\documentclass[11pt,letterpaper]{article}
\usepackage[margin=1.13in]{geometry}
\usepackage[T1]{fontenc}
\usepackage{lmodern,microtype,amsmath,amssymb,amsthm,mathtools}
\usepackage{graphicx,booktabs,xcolor,tikz}
\usetikzlibrary{arrows.meta}
\usepackage[colorlinks=true,linkcolor=teal!55!black,citecolor=teal!55!black,urlcolor=teal!55!black]{hyperref}
\hypersetup{pdftitle={When Do Sketches Achieve Gaussian Accuracy in Randomized SVD? A Geometric Perspective},pdfauthor={Yuning Yang}}
\newtheorem{theorem}{Theorem}[section]
\newtheorem{corollary}[theorem]{Corollary}
\newtheorem{lemma}[theorem]{Lemma}
\numberwithin{equation}{section}
\newcommand{\R}{\mathbb R}
\newcommand{\E}{\mathbb E}
\newcommand{\op}{\mathrm{op}}
\DeclareMathOperator{\tr}{tr}
\DeclareMathOperator{\Var}{Var}
\DeclareMathOperator{\diag}{diag}
\DeclareMathOperator{\range}{range}
\DeclareMathOperator{\Gr}{Gr}
\DeclareMathOperator{\SO}{SO}
\providecommand{\doi}[1]{\href{https://doi.org/#1}{\nolinkurl{#1}}}
\title{When Do Sketches Achieve Gaussian Accuracy in Randomized SVD? A Geometric Perspective}
\author{Yuning Yang\thanks{School of Mathematics, Guangxi University,
Nanning 530004, China. Email: \href{mailto:yyang@gxu.edu.cn}{yyang@gxu.edu.cn}.}}
\date{}
\begin{document}
\maketitle
\begin{abstract}
Fast sketches often achieve nearly Gaussian accuracy
 in randomized SVD experiments, despite having vastly
  different worst-case behavior. When and why can such
   agreement be expected? The mean approximation error
    of a sketch can vary across right singular bases, 
    and its average over these bases can differ from
     the Gaussian mean. For Gaussian sketches, both the
      variance and the bias are zero. Under Haar measure, 
      the bias also vanishes for general sketches, but 
      the variance need not. We establish non-asymptotic
      concentration bounds around the Gaussian mean using the
        Brascamp--Lieb inequality and give a
      geometric characterization of asymptotic agreement.
       At fixed target rank $k$ and oversampling $p\ge2$,
        the mean relative error approaches its Gaussian 
        counterpart in $L^1$ over Haar right singular 
        bases, uniformly over spectra, if and only if
         $\|\E P_{\mathcal{S}}\|_{\op}\to0$ as the input dimension 
         $n$ grows, where $P_{\mathcal{S}}$ is the orthogonal projector
          onto the sketch space. For $p\ge4$, the same characterization
          holds in $L^2$, and both the $L^1$ and $L^2$ discrepancies
          are bounded by $C_{k,p}\|\E P_{\mathcal{S}}\|_{\op}^{1/2}$. This gives
            an $O(n^{-1/2})$ bound for sketches such as subsampled randomized Walsh--Hadamard and cosine transforms and, more remarkably, for uniform coordinate sampling without any coordinate mixing.
The results also extend to a class of non-Haar input distributions.

\par\medskip\noindent
\textbf{Keywords:} randomized SVD, randomized sketching,  
mean error matching, Haar measure, Brascamp--Lieb inequality.
\end{abstract}

\section{Introduction}\label{sec:intro}

A recurring observation in randomized SVD experiments is that
fast sketches often achieve nearly the same accuracy as
Gaussian sketches. In an early
range-finder experiment, Halko, Martinsson, and Tropp even obtained a
slightly smaller spectral-norm error with subsampled Fourier sketches than
with Gaussian sketches \cite[Section~7.4, Figure~7.8]{halko2011}. The input
mattered, however. When Boutsidis and Gittens held the singular values
fixed and varied the right singular vectors, the Hadamard error fell on
the less coherent input, while the Gaussian error changed much less
\cite[Section~6.2]{boutsidis2013}.
Saibaba and Mi{\k{e}}dlar varied the decay of the singular values and also
found comparable errors across several sketch distributions on sea-surface
temperature data \cite[Section~5]{saibaba2025}. Taken together, these
experiments suggest that agreement can depend both on the decay of the
singular values and on the directions that carry them.

Martinsson and Tropp note that embeddings which mix coordinates often
behave much like Gaussian embeddings, while plain coordinate sampling
can perform markedly worse \cite[Warning~9.5]{martinssonTropp2020}.
Epperly and Webber varied both the spectrum and the right singular basis
in a single experiment
\cite[Section~3.2, Figure~2]{epperlyWebber2026}. They paired a step spectrum
and a polynomially decaying spectrum with either the coordinate or the
cosine basis for the right singular vectors. They tested
Gaussian, uniform, and sign matrices, two sparse sketches, and a
rerandomized cosine transform. By width $20$, all six gave nearly the same
mean squared Frobenius residual on every matrix. At smaller widths,
however, the sign and sparse sketches gave much larger errors on the step
spectrum in the coordinate basis.

\begin{figure}[tb]
\centering
\includegraphics[width=\linewidth]{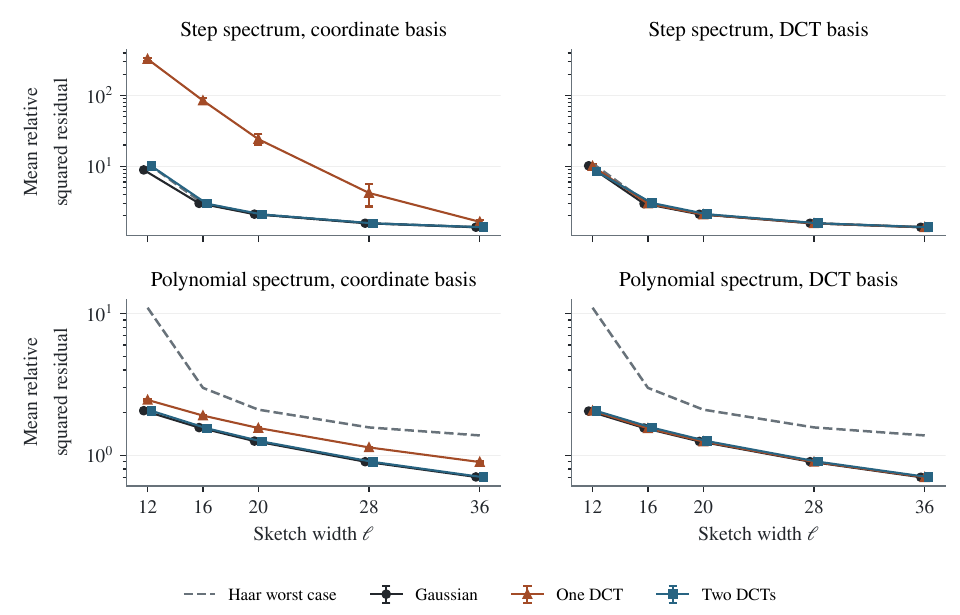}
\caption{One and two signed cosine passes compared with Gaussian sampling
on four prescribed inputs ($n=2000$, target rank $k=10$). The top row
has ten unit singular values and a flat tail of $10^{-3}$, and the bottom
row uses $\sigma_i=1/i$. The right singular basis is the coordinate
basis (left) or the DCT-II basis (right). Each point is the mean of 400
draws of the relative squared error in \eqref{eq:residual}. Bars show one
empirical standard error and the dashed line gives the Gaussian worst-case
mean in \eqref{eq:worst-case}. The draws are paired across panels and
widths.}\label{fig:opening}
\end{figure}

Figure~\ref{fig:opening} follows their four-panel design in a separate
experiment that singles out the effect of a second signed discrete
cosine transform (DCT). At width $20$ in the upper-left panel, the sample
mean squared residual for one cosine pass is $24.21$ times the squared
error of the best rank-ten approximation. Adding a second pass on the same
matrix brings this ratio down to $2.09$, close to the $2.11$ observed for
Gaussian sampling. The upper-right panel offers another clue. There the
singular values are exactly the same, yet changing the right singular
basis brings even one pass close to the Gaussian curve.

Does the close agreement observed in Figure~\ref{fig:opening}
extend to worst-case mean errors? To make this comparison precise,
let $A$ be an input with $n$ columns, and let $\mathcal{S} \coloneqq \range\Omega$
be the column space of the random sketch matrix $\Omega$. Write $P_{A\mathcal{S}}$ for
orthogonal projection onto the sampled range $A\mathcal{S}$. The range finder keeps
the whole sampled range \cite{halko2011,martinssonTropp2020}. We measure
its squared error relative to the tail energy,
\begin{equation}\label{eq:residual}
 \mathcal R_k(A,\mathcal{S}) \coloneqq \frac{\|(I-P_{A\mathcal{S}})A\|_F^2}{\tau_k^2},
 \qquad \tau_k^2 \coloneqq \sum_{i>k}\sigma_i^2>0,
\end{equation}
where $\sigma_1\ge\cdots\ge\sigma_n\ge0$ are the singular values, padded
with zeros if necessary.

Fix $1\le k<\ell<n$ and $p \coloneqq \ell-k\ge2$.
For the real Gaussian sketch space $\mathcal S_G$, Epperly and Webber
\cite[Theorems~3.2--3.3]{epperlyWebber2026}
established a sharp worst-case mean error formula.
In contrast, the finite support of the two-pass DCT sketch space
$\mathcal S_2$ leads to an unbounded worst-case mean error
(Appendix~\ref{app:paired-dct}). Specifically,
\begin{equation}\label{eq:worst-case}
 \sup_A\E\mathcal R_k(A,\mathcal{S}_G)
   =\frac{n-\ell}{n-k}\left(1+\frac{k}{p-1}\right),
 ~~{\rm while}~~
 \sup_A\E\mathcal R_k(A,\mathcal{S}_2)=\infty .
\end{equation}
Both suprema range over full-column-rank inputs. A
two-level spectrum approaches the Gaussian supremum as its spectral gap
grows \cite[Theorem~3.2 and Section~6]{epperlyWebber2026}.
To see why the second supremum is infinite, choose a space that the
two-pass DCT sketch draws with positive probability. An input can place a
leading right singular vector perpendicular to this space, so the sketch
misses that direction whenever it draws this space. Increasing the
corresponding singular value then makes the mean relative error arbitrarily
large while leaving the tail energy fixed.

{
The worst-case result tells us how large the mean error can be, but it
leaves open how often the close agreement in Figure~\ref{fig:opening}
occurs. In the upper panels, changing only the right singular basis
moves the one-pass DCT curve substantially, while the Gaussian mean is
the same for both bases. For a given spectrum and sketch, is agreement
with the Gaussian mean common or exceptional among inputs with those
singular values?
\par}

\subsection{\texorpdfstring{Mean error matching on a fixed input}{Mean error matching on a fixed input}}

{
In Figure~\ref{fig:opening}, each mean is estimated by redrawing the sketch
on the same input. Accordingly, fix the singular values
$\Sigma\coloneqq\diag(\sigma_1,\ldots,\sigma_n)$ with $\tau_k^2>0$ and
a right singular basis $V\in O(n)$; we may write the input as $A=\Sigma V^T$,
since the left singular vectors do not affect \eqref{eq:residual} and we omit them.
For a random sketch matrix $\Omega\in\mathbb R^{n\times\ell}$ with full
column rank almost surely, let $\mathcal S\coloneqq\range\Omega$ denote
its column space. Particularly  for a standard Gaussian
matrix $G$, write $\mathcal S_G\coloneqq\range G$. Since even independent
draws from the same sketch distribution may produce different errors on
a fixed input, we compare different sketches through their mean
approximation errors. On this input, the two mean errors are
\begin{equation}\label{eq:means}
\begin{aligned}
 F_\Sigma(V)&\coloneqq\E_{\mathcal S}\mathcal R_k(\Sigma V^T,\mathcal S),\\
 g_\Sigma&\coloneqq\E_G\mathcal R_k(\Sigma V^T,\mathcal S_G)
 =\E_G\mathcal R_k(\Sigma,V^T\mathcal S_G)
 =\E_G\mathcal R_k(\Sigma,\mathcal S_G),
\end{aligned}
\end{equation}
where the last equality follows from rotational invariance, so $g_\Sigma$
is the same for every right singular basis $V$.
Replacing the Gaussian sketch by $\mathcal S$ then changes the error,
on average, by
\[
 \E_{\mathcal S,G}\!\left[
   \mathcal R_k(A,\mathcal S)-\mathcal R_k(A,\mathcal S_G)
 \right]=F_\Sigma(V)-g_\Sigma .
\]
In contrast to the Gaussian mean $g_\Sigma$, the mean error $F_\Sigma(V)$ of a general sketch may depend on the right singular basis $V$.

Agreement on a few prescribed inputs does not tell us how common
such agreement is among inputs with the same spectrum. We therefore allow the right singular basis $V$ to follow
a probability distribution independent of the sketch. We are interested
in the probability that the two means agree to within a given tolerance
$\varepsilon>0$:
\[
 \Pr_V\!\left\{|F_\Sigma(V)-g_\Sigma|\le\varepsilon\right\}.
\]
For any such distribution, the bias--variance decomposition gives
\begin{equation}\label{eq:mean-discrepancy}
 \E_V|F_\Sigma(V)-g_\Sigma|^2
 =\Var_V F_\Sigma(V)+|\E_V F_\Sigma(V)-g_\Sigma|^2 .
\end{equation}
The first term measures variation of the sketch mean across bases,
and the second is its squared bias relative to the Gaussian mean.
For a Gaussian sketch, both terms vanish regardless of the distribution
of $V$. For a general sketch, however, neither term need vanish. By Markov's
inequality,
\[
 \Pr_V\{|F_\Sigma(V)-g_\Sigma|>\varepsilon\}
 \le \frac{\Var_V F_\Sigma(V)+|\E_V F_\Sigma(V)-g_\Sigma|^2}
          {\varepsilon^2}.
\]
Thus controlling both terms ensures agreement on most inputs under the
chosen distribution.

The distribution of the right singular basis may depend on the
application. As a starting point, we consider the model
that does not favor any particular orientation of the right singular
vectors. This leads to the uniform distribution of $V$,
namely Haar probability  on $O(n)$.
Nevertheless, our results also extend to a class of non-Haar distributions.
For Haar distribution, we draw $V$  once and hold it fixed across sketch
draws. For every sketch law, rotational invariance gives
\[
 \begin{aligned}
 V^T\mathcal S \overset{d}{=}\mathcal S_G
 \Longrightarrow
 \mathcal R_k(\Sigma V^T,\mathcal S)
 \overset{d}{=}\mathcal R_k(\Sigma,\mathcal S_G)
 \Longrightarrow 
 \E_V F_\Sigma(V) =g_\Sigma .
 \end{aligned}
\]
\cite[Section~6.4]{epperlyWebber2026} uses such a rotation-averaging argument in their
proof of the minimax lower bound.
However, mean error matching on individual inputs remains a nontrivial
problem even under the Haar input model: The identity
$\E_V F_\Sigma(V)=g_\Sigma$ holds for every sketch law, but does not
guarantee that $F_\Sigma(V)$ is close to $g_\Sigma$ for most choices of $V$.
In fact, under the Haar model, the second term of \eqref{eq:mean-discrepancy} vanishes, but the first may not, reducing \eqref{eq:mean-discrepancy} to
\[
 \E_V|F_\Sigma(V)-g_\Sigma|^2=\Var_V F_\Sigma(V).
\]
Therefore, mean error matching becomes a concentration problem for
$F_\Sigma(V)$ around $g_\Sigma$.

Since $V^T\mathcal S\overset{d}{=}\mathcal S_G$, we have 
$
 \Var_G\mathcal R_k(\Sigma,\mathcal S_G)
 =
 \Var_{V,\mathcal S}\mathcal R_k(\Sigma V^T,\mathcal S)
$.
Applying the law of total variance to the right-hand side then gives
\begin{equation}\label{eq:total-variance}
 \Var_G\mathcal R_k(\Sigma,\mathcal S_G)
 =
 \Var_{V,\mathcal S}\mathcal R_k(\Sigma V^T,\mathcal S)=\E_V\Var_{\mathcal S}\mathcal R_k(\Sigma V^T,\mathcal S)
 +\Var_V F_\Sigma(V).
\end{equation}
This shows that  under the Haar   model, the left-hand side of (1.5) is independent of the sketch law, whereas the two terms on the right-hand side can vary substantially. To see this, consider two extreme cases:
\[
 \begin{aligned}
 \text{Gaussian sketch:}\qquad
 &\Var_V F_\Sigma(V)=0,\\
 \text{Fixed first $\ell$ coordinates:}\qquad
 &\E_V\Var_{\mathcal S}\mathcal R_k(\Sigma V^T,\mathcal S)=0,\\
 &\Var_V F_\Sigma(V)=\Var_G\mathcal R_k(\Sigma,\mathcal S_G).
 \end{aligned}
\]
In the first case, all variance comes from sketch draws; in the second,
all comes from changing the input basis. Moreover, for the second
case, the variance of $F_\Sigma(V)$ need not vanish even as $n$ grows
at fixed sketch width (Section~\ref{sec:projection-limits}).
As a result, the above shows that the concentration of the sketch mean
depends crucially on the sketch law, even under the Haar   model.

In summary, our aim is to characterize which sketch laws match the Gaussian mean on typical inputs and to quantify the difference in finite dimensions, uniformly over   spectra.
\par}

\subsection{A geometric characterization of Gaussian mean error matching}

{
We work under the Haar model of the preceding subsection.
Fix the target rank $k$ and oversampling $p\ge2$, so that the sketch
width $\ell=k+p$ stays fixed as $n$ grows.
To characterize this concentration geometrically, we measure how strongly
the sketch law can favor a particular direction and compare the mean errors
uniformly over spectra.
Let $P_{\mathcal{S}}$ be the orthogonal projector onto the sketch space $\mathcal{S}$,
and define
\begin{equation}\label{eq:rho}
 \rho_n \coloneqq \|\E P_{\mathcal{S}}\|_{\op}
       =\max_{\|u\|=1}\E\|P_{\mathcal{S}}u\|^2,\qquad
 D_n \coloneqq \sup_{\Sigma:\,\tau_k>0}\E_V\,|F_\Sigma(V)-g_\Sigma|.
\end{equation}
$\rho_n$ depends only on the sketch law and gives the maximum
expected energy captured from a unit vector. For fixed sampling of the
first $\ell$ coordinates, $\rho_n=1$, whereas for uniform coordinate
sampling, $\rho_n=\ell/n$.
\par}
\begin{theorem}\label{thm:criterion}
{
Fix \(k\ge1\), \(p\ge2\), and \(\ell \coloneqq k+p\). For each \(n\ge\ell\), let
\(\mathcal{S}\) be an arbitrarily distributed random \(\ell\)-dimensional subspace
of \(\R^n\), independent of the Haar   singular basis \(V\).

\smallskip\noindent
\textup{(i) Variance comparison.}
For every spectrum with \(\tau_k>0\),
\begin{equation}\label{eq:variance-comparison}
 \Var_V F_\Sigma(V)
 \le 2\rho_n\,\Var_G\mathcal R_k(\Sigma,\mathcal{S}_G).
\end{equation}

\smallskip\noindent
\textup{(ii) Uniform \(L^1\) agreement.}
For any $r$ satisfying \(1<r\le2\) and \(2r<p+1\), one has
\(D_n\le C_{k,p,r}\rho_n^{\,1-1/r}\), with $C_{k,p,r}$ independent of \(n\)
and the spectrum. In particular, for \(p\ge4\),
\[
 D_n\le C_{k,p}\sqrt{\rho_n}.
\]
Moreover, as \(n\to\infty\) with \(k,p\) fixed,
\begin{equation}\label{eq:criterion}
             \rho_n\longrightarrow0
       \quad\Longleftrightarrow\quad D_n\longrightarrow0 .
\end{equation}

\smallskip\noindent
\textup{(iii) Uniform \(L^2\) agreement.}
For \(p\ge4\),
\begin{equation}\label{eq:uniform-second-moment}
 \sup_{\Sigma:\,\tau_k>0}\E_V|F_\Sigma(V)-g_\Sigma|^2
 \le C_{k,p}^2\rho_n.
\end{equation}
Moreover, as \(n\to\infty\) with \(k,p\) fixed,
\[
 \rho_n\longrightarrow0
 \quad\Longleftrightarrow\quad
 \sup_{\Sigma:\,\tau_k>0}\E_V|F_\Sigma(V)-g_\Sigma|^2
 \longrightarrow0 .
\]
\par}
\end{theorem}

{
\paragraph{Discussions.} The non-asymptotic   comparison \eqref{eq:variance-comparison} together with
\eqref{eq:total-variance}  bounds the relative contribution of
changes in the input basis to the total error variance:
$
 \frac{\Var_V F_\Sigma(V)}{\Var_G\mathcal R_k(\Sigma,\mathcal S_G)}
 \le \min\{1,2\rho_n\}
$. 
Thus, as $\rho_n\to0$, the second term in (1.5) becomes negligible relative to the total variance, which is then dominated by sketch randomness.
For $p\ge4$, the uniform second-moment bound gives, by Chebyshev's
inequality,
\[
 \sup_{\Sigma:\,\tau_k>0}
 \Pr_V\{|F_\Sigma(V)-g_\Sigma|>\varepsilon\}
 \le C_{k,p}^2\rho_n/\varepsilon^2 .
\]
For every $p\ge2$, Markov's inequality also bounds this probability
by $D_n/\varepsilon$. These bounds are uniform over spectra, although
the exceptional set of bases may depend on the spectrum.
\par}

To establish the necessity of $\rho_n\to0$ for $D_n\to0$ in the theorem,
it suffices to consider the spectrum
$\Sigma=\diag(\sqrt{n+1},1,\ldots,1)$.
The proof also yields the following lower bound for $D_n$.
Together with Theorem~\ref{thm:criterion}, we have, for all
sufficiently large $n$,
\begin{equation}\label{eq:two-sided}
 c_\ell\left(\rho_n-{\ell}/{n}\right)^2\le D_n
 \qquad (p\ge2),\qquad
 D_n\le C_{k,p}\sqrt{\rho_n}\qquad(p\ge4).
\end{equation}
The gap between the upper and lower bounds leaves open
the optimal convergence rate.

The key tool for establishing the variance bound in Theorem~\ref{thm:criterion} 
is the rotation-group Brascamp--Lieb inequality of Barthe et al. \cite{barthe2011}, whose geometric condition naturally leads to the quantity $\rho_n$. This inequality will be introduced in Section~\ref{sec:prelims}.

{
 \paragraph{Deterministic $V$.} Can the geometric criterion in Theorem~\ref{thm:criterion} guarantee
mean error matching for arbitrary deterministic $V$? The following
example shows that this is not the case. For any prescribed
 $V_0\in O(n)$, choose $\mathcal{S}$ as the span of
$\ell$ columns of $V_0$ sampled uniformly without replacement. Then
\[
 \E P_{\mathcal{S}}=\frac{\ell}{n}V_0V_0^T=\frac{\ell}{n}I_n,
\]
so this sketch attains the smallest possible value of $\rho_n$.
Nevertheless, for $A=\Sigma V_0^T$, we have
$
 \|(I-P_{A\mathcal S})A\|_F^2
 =\sum_{j:\,V_0e_j\notin\mathcal S}\sigma_j^2
$. 
Since $\Pr(V_0e_j\notin\mathcal S)=1-\ell/n$ for every $j$,
$
 F_\Sigma(V_0)
 =\frac{1-\ell/n}{\tau_k^2}\sum_{j=1}^n\sigma_j^2
$. 
Now taking $k=1$ and
$\Sigma=\diag(\sqrt{\Gamma},1,\ldots,1)$ with $\Gamma>1$
gives
\[
 F_\Sigma(V_0)
 =\frac{n-\ell}{n-1}
  \left(1+\frac{\Gamma-1}{n}\right).
\]
Letting $\Gamma=n^2$, this mean diverges as $n$ grows at fixed sketch
width, whereas the Gaussian mean remains bounded by
\eqref{eq:worst-case}. Thus for $V\equiv V_0$, the first term on the right-hand side of
\eqref{eq:mean-discrepancy} is zero, whereas the second term diverges.
Consequently, $\E_V|F_\Sigma(V)-g_\Sigma|^2$ also diverges. 

\par}

{
\paragraph{Non-Haar distributions.} We can extend
Theorem~\ref{thm:criterion} to a class of non-Haar distributions.
Let $\nu_n$ denote the Haar   measure on $O(n)$ and
$\lambda_n$ be another probability measure on $O(n)$. We write
$\lambda_n\le C_n\nu_n$ if
\[
 \lambda_n(B)\le C_n\nu_n(B)
 \qquad\text{for every Borel set } B\subseteq O(n).
\]
\begin{corollary}\label{cor:non-haar}
Under the other hypotheses of Theorem~\ref{thm:criterion}, let
$V\sim\lambda_n$ be independent of the sketch, with
$\lambda_n\le C_n\nu_n$. Then
\begin{equation}\label{eq:non-haar-bounds}
 \begin{aligned}
 \sup_{\Sigma:\,\tau_k>0}\E_{V\sim\lambda_n}|F_\Sigma(V)-g_\Sigma|
 &\le C_n D_n,\\
 \E_{V\sim\lambda_n}|F_\Sigma(V)-g_\Sigma|^2
 &\le 2C_n\rho_n\,\Var_G\mathcal R_k(\Sigma,\mathcal S_G)
 \quad (\tau_k>0),
 \end{aligned}
\end{equation}
where $D_n$ is still defined with respect to $\nu_n$ in
\eqref{eq:rho}. If $C_n$ is bounded independently of $n$, the
sufficient direction and its convergence rates are preserved.
Furthermore, if $c\nu_n\le\lambda_n\le C\nu_n$ for constants $0<c\le C<\infty$
independent of $n$, then the full characterization holds:
\begin{equation}\label{eq:non-haar-criterion}
 \rho_n\to0
 \quad\Longleftrightarrow\quad
 \sup_{\Sigma:\,\tau_k>0}\E_{V\sim\lambda_n}
       |F_\Sigma(V)-g_\Sigma|\to0 .
\end{equation}
\end{corollary}

As in Theorem \ref{thm:criterion}   $L^2$ agreement can also be obtained. The non-Haar distributions in Section~\ref{sec:non-haar} illustrate
that this extension permits genuine directional preferences.
\par}

\subsection{Implications from Theorem~\ref{thm:criterion}}

\noindent\textbf{Gaussian mean error for orthogonal transforms.}
Uniformly sampling rows of an orthogonal matrix gives
\[
  \E P_{\mathcal{S}}=\frac{\ell}{n}I_n \qquad\Longrightarrow\qquad D_n\to0 .
\]
Thus uniform coordinate sampling, the subsampled randomized Hadamard
transform, and one or more signed cosine passes all match the Gaussian
mean asymptotically (Corollary~\ref{cor:orthogonal-rows}).

{
For these sketches, Theorem~\ref{thm:criterion} further gives
$D_n=O(n^{-1/2})$ for $p\ge4$ at fixed sketch width.
\par}

\medskip
\noindent\textbf{Gaussian mean error beyond orthogonal transforms.}
Under suitable symmetry conditions, non-Gaussian sketches have the same
mean projector,
\[
  \E P_{\mathcal{S}}=\frac{\ell}{n}I_n .
\]
{
The same conclusions hold for sketches with iid uniform entries and for
the sign and sparse sketches in Corollary~\ref{cor:symmetries},
conditioned on full column rank in the discrete cases.
\par}

\medskip
\noindent\textbf{Simultaneous agreement on the same input.}
For a fixed collection of $m$ sketches satisfying the criterion, their
conditional means $F_\Sigma^{(j)}(V)$ obey
\[
  \sup_{\Sigma:\,\tau_k>0}\E_V\max_{1\le j\le m}
  |F_\Sigma^{(j)}(V)-g_\Sigma|\longrightarrow0 .
\]
Thus even the largest of the $m$ absolute differences is small on a typical
input shared by all the sketches (Corollary~\ref{cor:common-mean}).

{
These examples show that a sketch need not resemble a Gaussian matrix
to match its mean error on typical inputs.
\par}

\subsection{Relation to previous results}

For a prescribed matrix, classical analyses of fast range finders control
how well a sampled transform embeds its leading right singular subspace
\cite{woolfe2008,tropp2011srht,boutsidis2013}. Saibaba and Mi{\k{e}}dlar
bound the error for several classes of sketches with independent columns
\cite{saibaba2025}. Cama{\~n}o, Epperly, Meyer, and Tropp give fixed-input
error guarantees for sparse sketches and study how many nonzeros they need
to compete with Gaussian sampling \cite[Theorem~1.3 and Figure~5]{camanoOSI}.

A different line of research investigates the distributional properties
of structured random transforms in comparison with Gaussian or Haar laws
\cite{haimiCohenLai,ahfock2021,zilcaMendelson}.
{
However, asymptotic Gaussianity of sketch projections does not guarantee
agreement in mean RSVD errors. As shown in Section~\ref{sec:discussion},
even a Gaussian projection limit can coexist with a diverging mean
relative error.
\par}

Universality has also been established for other problems in randomized
numerical linear algebra. Oymak and Tropp derive universality laws for
the success threshold and stability of randomized dimension reduction
\cite{oymakTropp}. In least-squares regression, Dobriban and Liu derive
common asymptotic variance and prediction efficiencies for Gaussian
and other iid sketches, and a separate pair of limits shared by Haar
projection and the subsampled randomized Hadamard transform. For
rotationally invariant designs, uniform sampling also attains the Haar
limits \cite{dobribanLiu}.

More closely related to our work are studies of the mean approximation error in randomized SVD. For Gaussian sketches, Akhtiamov, Ghane, and Hassibi derive an asymptotically exact error formula in terms of the full singular value spectrum as the sketch width grows \cite[Definition~1 and Theorem~2]{akhtiamovRSVD}. Under a stable-rank condition, Derezi\'nski, Liang, Liao, and Mahoney derive a common approximation for the mean squared Frobenius error of Gaussian and other iid sub-Gaussian sketches, including random sign matrices \cite[Theorem~1 and Corollary~1]{derezinskiPrecise}. Building on this result, Derezi\'nski shows that sparse LESS sketches can nearly match the mean error of Gaussian sketches of the same width \cite[Corollary~15]{derezinskiLESS}.

These results provide theoretical explanations for the agreement between particular non-Gaussian sketches and Gaussian sketches. However, a natural question remains: which sketch laws exhibit such
agreement on typical inputs, uniformly over singular value spectra?
Theorem~\ref{thm:criterion} provides a necessary and sufficient
geometric condition for such agreement in the Haar   model, together
with non-asymptotic bounds on the difference between the two means.

The remainder is organized as follows: Section~\ref{sec:prelims} explains the geometry of the mean projector and
the Brascamp--Lieb inequality for a shared rotation, which Section~\ref{sec:proof} uses
to prove the criterion and its extension to non-Haar
distributions. Section~\ref{sec:applications} derives  
corollaries for common sketches and illustrates them numerically, while
Section~\ref{sec:discussion} examines what can happen when a singular
direction is chosen in advance and explores the scope of
the non-Haar extension.

\section{Geometry and preliminaries}\label{sec:prelims}

\subsection{What the mean projection measures}

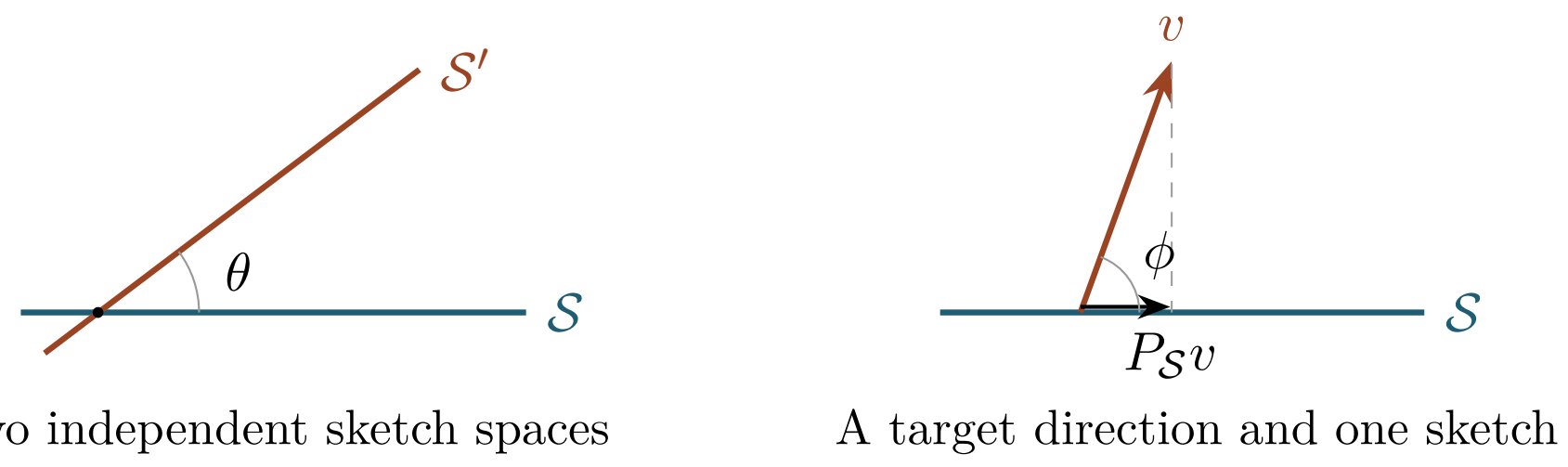
\begin{figure}[!htbp]
\centering
\begin{tikzpicture}[x=1cm,y=1cm,>=Stealth]
 \definecolor{slate}{RGB}{31,93,116}\definecolor{rust}{RGB}{155,68,36}
 \begin{scope}[shift={(0,0)}]
 \draw[slate,very thick] (-.55,0)--(3.05,0) node[right] {$\mathcal{S}$};
 \draw[rust,very thick] (-.38,-.29)--(2.29,1.73) node[right] {$\mathcal{S}'$};
 \draw[gray!80] (.72,0) arc[start angle=0,end angle=37,radius=.72];
 \node at (1.00,.29) {$\theta$};
 \fill (0,0) circle (1.1pt);
 \node[below] at (1.26,-.56) {\small Two independent sketch spaces};
 \end{scope}
 \begin{scope}[shift={(7,0)}]
 \draw[slate,very thick] (-1,0)--(2.45,0) node[right] {$\mathcal{S}$};
 \draw[rust,very thick,->] (0,0)--(.65,1.79) node[above] {$v$};
 \draw[gray!75,dashed] (.65,1.77)--(.65,0);
 \draw[black,thick,->] (0,.04)--(.64,.04);
 \node[below] at (.64,-.01) {$P_{\mathcal{S}}v$};
 \draw[gray!80] (.42,0) arc[start angle=0,end angle=70,radius=.42];
 \node at (.56,.41) {$\phi$};
 \node[below] at (.83,-.56) {\small A target direction and one sketch};
 \end{scope}
\end{tikzpicture}
\caption{Two geometric quantities for one-dimensional sketch spaces:
$\cos^2\theta$ measures the overlap between two sketch spaces (left),
and $\cos^2\phi=\|P_{\mathcal S}v\|^2$ measures the captured energy
of a unit vector $v$ (right).}
\label{fig:angles}
\end{figure}

{
For a unit vector \(u\), the expected energy captured by the sketch is
$
 \E\|P_{\mathcal S}u\|^2=u^TMu$, with $M\coloneqq\E P_{\mathcal S}
$.
Since \(0\preceq M\preceq I\) and \(\tr M=\ell\), one has
\[
 \frac{\ell}{n}=\frac{\tr M}{n}\le\rho_n=\lambda_{\max}(M)\le1.
\]
Let \(\mathcal S,\mathcal S'\) be independent copies of the sketch space,
with orthonormal basis matrices \(Q,Q'\), respectively. Their principal
angles \(\theta_1,\ldots,\theta_\ell\) satisfy
\(\cos\theta_j=\sigma_j(Q^TQ')\) \cite{bjorckGolub1973}, so
\[
 \sum_{j=1}^{\ell}\cos^2\theta_j
 =\|Q^TQ'\|_F^2=\tr(P_{\mathcal S}P_{\mathcal S'}).
\]
By independence and the bounds on the eigenvalues of \(M\),
\begin{equation}\label{eq:overlap}
 \E\sum_{j=1}^{\ell}\cos^2\theta_j
 =\tr(M^2),\qquad
 \rho_n^2\le\tr(M^2)\le\ell\rho_n.
\end{equation}
Consequently, at fixed sketch width,
\[
 \rho_n\to0
 \quad\Longleftrightarrow\quad
 \E\sum_{j=1}^{\ell}\cos^2\theta_j\to0.
\]
Figure~\ref{fig:angles} distinguishes the angle between two sketch spaces
from that between a sketch space and a target direction.
\par}

{
Related criteria appear in Gaussian projection theory
\cite{diaconisFreedman}. D{\"u}mbgen and Del Conte-Zerial
characterize standard Gaussian projection limits through concentration of
normalized lengths and vanishing normalized inner products of independent
draws \cite[Corollary~2.2]{dumbgenDelConte}. Here we instead study
sketch-space overlap and normalized residuals.

For coordinate sampling, as  \(P_{\mathcal S}\) is diagonal, and then
$
 \rho_n=\max_{1\le j\le n}\Pr(e_j\in\mathcal S)
$.
The criterion can hold even when only \(o(n)\) coordinates have positive
inclusion probability. For example, sample \(\ell\) distinct coordinates
uniformly from the first \(m=\lfloor\sqrt n\rfloor\ge\ell\). Then
\[
 \rho_n=\frac{\ell}{m}\to0,\qquad \frac{m}{n}\to0,\qquad
 P_{\mathcal S}e_j=0\quad(j>m).
\]
\par}

\subsection{The error as a function of a subspace}

{
Let \(\Gr(\ell,n)\) denote the Grassmannian of \(\ell\)-dimensional
subspaces of \(\R^n\). Since right multiplication by \(V\) preserves
the Frobenius norm, for \(A=\Sigma V^T\) we have
\begin{equation}\label{eq:h}
 \mathcal R_k(\Sigma V^T,\mathcal S)
 =\mathcal R_k(\Sigma,V^T\mathcal S).
\end{equation}
The range of a standard Gaussian matrix is uniform on \(\Gr(\ell,n)\)
by orthogonal invariance and almost sure full column rank \cite{mezzadri}.
Write \(\mu\) for the law of \(\mathcal S\) and \(\nu\) for the uniform law
on \(\Gr(\ell,n)\). For each fixed \(\mathcal S\), \(V^T\mathcal S\) has law
\(\nu\) when \(V\) is Haar. Thus \eqref{eq:h} gives
\begin{equation}\label{eq:center}
 \begin{aligned}
 F_\Sigma(V) =\int\mathcal R_k(\Sigma,V^T\mathcal S)\,d\mu(\mathcal S),~~
 \E_V F_\Sigma(V)&=\int\mathcal R_k(\Sigma,\mathcal S)\,d\nu(\mathcal S)
 =g_\Sigma .
 \end{aligned}
\end{equation}
The bound \(0\le\mathcal R_k(\Sigma,\mathcal S)\le
\|\Sigma\|_F^2/\tau_k^2\) justifies the interchange of integrals for each
spectrum. However, this bound is not uniform over spectra, so we will
need a uniform moment estimate. Zero singular values are allowed:
\(P_{\Sigma\mathcal S}\) is a measurable function of \(\mathcal S\), since
\[
 P_{\Sigma\mathcal S}=\lim_{\varepsilon\downarrow0}
 \Sigma P_{\mathcal S}\Sigma
 (\Sigma P_{\mathcal S}\Sigma+\varepsilon I)^{-1}.
\]
\par}

For the converse, we will also compute a residual in the original
coordinates. Let \(K \coloneqq A^TA\), and let \(Q\) be an orthonormal
basis matrix for \(\mathcal S=\range\Omega\), i.e.,
\(Q^TQ=I_\ell\) and \(\range Q=\mathcal S\).
Then \(P_{A\mathcal S}=AQ(Q^TKQ)^\dagger Q^TA^T\),
where \({}^\dagger\) denotes the
Moore--Penrose inverse. The total squared energy is \(\tr K\).
Cyclically permuting factors in the trace of the captured energy gives
\begin{equation}\label{eq:captured}
 \|P_{A\mathcal S}A\|_F^2=\tr\!\big[(Q^TKQ)^\dagger Q^TK^2Q\big].
\end{equation}
This is the trace form of the Gram-matrix correspondence used by Epperly and
Webber \cite[Sections~3.4 and 6.1]{epperlyWebber2026}.

\subsection{The Brascamp--Lieb inequality for a shared rotation}\label{sec:rotation-tool}

Imagine rotating two orthogonal lines in the plane by the same random
angle. Each line separately is uniform, yet observing the first determines
the second. Rotating a family of sketch spaces likewise preserves all the
angles between them.

{
\paragraph{Brascamp--Lieb inequality.}
Let \(\mathcal S_1,\ldots,\mathcal S_m\in\Gr(\ell,n)\), let \(c_i>0\)
for \(i=1,\ldots,m\), and let \(U\) be Haar distributed on
\(\SO(n)\), the determinant-one orthogonal group.
Suppose that
\begin{equation}\label{eq:directional-energy}
 \sum_i\nolimits c_i\|P_{\mathcal S_i}u\|^2\le\|u\|^2\quad\text{for all }u\in\R^n
 \qquad\Longleftrightarrow\qquad
 \sum_i\nolimits c_iP_{\mathcal S_i}\preceq I .
\end{equation}
Then the rotation-group Brascamp--Lieb inequality of Barthe et al.
\cite[Theorem~7, Eq.~(13), and Proposition~10]{barthe2011} gives
\begin{equation}\label{eq:BCLM}
 \E_U\prod_i\nolimits H_i(U\mathcal S_i)^{c_i/2}
 \le \prod_i\nolimits \big(\E_U H_i(U\mathcal S_i)\big)^{c_i/2}
\end{equation}
for any nonnegative bounded Borel functions
\(H_i:\Gr(\ell,n)\to[0,\infty)\).

For mutually orthogonal sketch spaces, \eqref{eq:directional-energy}
holds with $c_i=1$; suitable weights also allow overlap.
This condition on the unrotated spaces brings the mean projector
in \eqref{eq:rho} into the proof.
To see how \eqref{eq:BCLM} can be used to bound a variance, let \(f\) be a
bounded real Borel function on \(\Gr(\ell,n)\) and take \(H_i=e^{2tf}\).
Then
\[
 \prod_i\nolimits H_i(U\mathcal S_i)^{c_i/2}
 =\exp\left(t\sum_i\nolimits c_i f(U\mathcal S_i)\right).
\]
Since, for any bounded real random variable \(Z\), $(\log\E e^{tZ})''|_{t=0}=\Var(Z)$,
applying \eqref{eq:BCLM} and comparing the second-order terms gives
a variance bound, as shown in Lemma~\ref{lem:variance}.

\paragraph{Derivation of \eqref{eq:BCLM}.}
We explain why the geometric condition \eqref{eq:directional-energy}
allows the exponents \(c_i/2\) in the   inequality.
The function \(H_i(U\mathcal S_i)\) is unchanged if \(U\) is replaced
by \(UR\), where \(R\in\SO(n)\) fixes \(\mathcal S_i\) pointwise.
The Lie algebra of this subgroup consists of skew-symmetric matrices
\(X\) satisfying \(XP_{\mathcal S_i}=0\). Within the space of
skew-symmetric matrices, the orthogonal projector (with respect to
the Frobenius inner product) onto the complement of this Lie algebra is
\[
 \Pi_i X=X-(I-P_{\mathcal S_i})X(I-P_{\mathcal S_i}).
\]
The rotation-group theorem admits exponents \(d_i\) whenever
\[
 \sum_i d_i\|\Pi_i X\|_F^2\le\|X\|_F^2
 \qquad\text{for every }X^T=-X.
\]
Using the block decomposition along \(\mathcal S_i\) and
\(\mathcal S_i^\perp\), one obtains
$
 \|\Pi_i X\|_F^2=2\|P_{\mathcal S_i}X\|_F^2
 -\|P_{\mathcal S_i}XP_{\mathcal S_i}\|_F^2
 \le 2\|P_{\mathcal S_i}X\|_F^2 
$. 
Consequently, \eqref{eq:directional-energy} gives
\[
 \begin{aligned}
 \sum_i\nolimits\frac{c_i}{2}\|\Pi_i X\|_F^2
 &\le\sum_i\nolimits c_i\|P_{\mathcal S_i}X\|_F^2 
  =\tr\left(X^T\left(\sum_i\nolimits c_iP_{\mathcal S_i}\right)X\right)
 \le\|X\|_F^2.
 \end{aligned}
\]
Thus   \(d_i=c_i/2\) gives \eqref{eq:BCLM} by the
rotation-group theorem \cite[Lemma~8 and Proposition~10]{barthe2011}.
\par}

\subsection{Two Gaussian calculations}\label{sec:gaussian-tools}

Let $B$ be the first $k$ rows of an $n\times\ell$ Gaussian sketch
and set $W\coloneqq BB^T$. For the $j$th row $b_j^T$ of $B$,
let $L_j$ be the span of the remaining rows in $\R^\ell$.
The Schur complement formula gives
\[
    \frac{1}{(W^{-1})_{jj}}=\|(I-P_{L_j})b_j\|^2 .
\]
Since $b_j\sim N(0,I_\ell)$ is independent of $L_j$ and
$\dim L_j^\perp=\ell-k+1=p+1$ almost surely, the right-hand side
has law $\chi^2_{p+1}$. For $Y$ of this law and $r>0$,
direct integration gives
\[
 \E Y^{-r}
 =2^{-r}\frac{\Gamma((p+1)/2-r)}{\Gamma((p+1)/2)}
 \quad\text{when }2r<p+1 .
\]
Applying convexity to the average of the \(k\) diagonal entries yields
\begin{equation}\label{eq:iw}
 \E(\tr W^{-1})^r\le k^r\,\E Y^{-r}
 \qquad(1\le r<(p+1)/2).
\end{equation}
Inverse Wishart moment formulas are given, for example, by Matsumoto
\cite{matsumoto2012}. Here the projection onto \(L_j^\perp\) explains
directly why \(p+1\) enters the exponent.

We also use a uniform unit vector \(v\) in \(\R^n\). Writing \(v \coloneqq g/\|g\|\)
for a standard Gaussian vector gives, for any fixed rank-\(d\) projector
\(P\), \(0<d<n\),
\begin{equation}\label{eq:sphere}
 v^TPv\ \stackrel{\mathrm{law}}{=}\  \frac{\sum_{j=1}^{d}g_j^2}{\sum_{j=1}^n g_j^2}
 \ \sim\ \operatorname{Beta}\bigl(d/2,(n-d)/2\bigr).
\end{equation}
The numerator and the remaining sum are independent chi-square variables.
In particular, \(d=1\) gives
\(\E(u^Tv)^2=1/n\) and \(\E(u^Tv)^4=3/[n(n+2)]\) for every unit \(u\).
If \(P\) is an independent random rank-\(\ell\) projector, the conditional
law of \(nv^TPv\) given \(P\) is the same for every \(P\). Thus \(nv^TPv\)
is independent of \(P\), which will let us separate the two averages
in \eqref{eq:spike-cov}.

\section{Proof of the criterion}\label{sec:proof}

We prove the three parts of Theorem~\ref{thm:criterion} in order.
First we establish the variance comparison. With a uniform moment
bound, we prove the \(L^1\) bounds and use a single-spike spectrum
for the converse. Finally, we establish the \(L^2\) conclusions.
Assume without loss of generality that $n>\ell$.

\subsection{Variance comparison}

We first prove part (i) of Theorem~\ref{thm:criterion} by establishing
a variance estimate for a shared rotation.

\begin{lemma}\label{lem:variance}
Let \(U\) be Haar on \(\SO(n)\) and \(f\) be a bounded real Borel function
on \(\Gr(\ell,n)\). Recall that \(\mathcal S_G=\range G\) has law \(\nu\).
For any sketch law \(\mu\),
\begin{equation}\label{eq:variance}
 \Var_U\!\left(\int f(U\mathcal{S})\,d\mu(\mathcal{S})\right)
 \le 2\rho_n\,\Var_G(f(\mathcal S_G)).
\end{equation}
\end{lemma}

\begin{proof}
Suppose first that the sketch takes just the spaces \(\mathcal{S}_1,\ldots,\mathcal{S}_m\),
with positive probabilities \(w_i\). Our aim is to recover the variance of
\(\sum_iw_i f(U\mathcal{S}_i)\) as a second derivative. To obtain these probabilities
as weights in \eqref{eq:BCLM}, let \(c_i \coloneqq w_i/\rho_n\). This choice is
admissible because \(\sum_i c_iP_{\mathcal{S}_i}=M/\rho_n\preceq I\), and
\(\rho_n\ge\ell/n>0\). Taking \(H_i\coloneqq e^{2\rho_n t f}\) now gives
\[
 \log\E_U \exp\!\left(t\sum_i\nolimits w_i f(U\mathcal{S}_i)\right)
 \le \frac{1}{2\rho_n}\log\E_G\exp(2\rho_n t f(\mathcal S_G)).
\]
Here every \(U\mathcal{S}_i\) has law \(\nu\), so the factors on the right of
\eqref{eq:BCLM} have a common expectation, and \(\sum_i c_i=1/\rho_n\).
To read a variance from this
inequality, recall that for a bounded real \(X\),
\[
       \log\E e^{tX}=t\E X+\tfrac12t^2\Var(X)+o(t^2).
\]
{
Since \(\E_U f(U\mathcal S_i)=\E_G f(\mathcal S_G)\)
and \(\sum_i w_i=1\), the linear terms on both sides agree.
Comparing the quadratic terms gives
\[
 \Var_U\left(\sum_i\nolimits w_i f(U\mathcal S_i)\right)
 \le 2\rho_n\Var_G(f(\mathcal S_G)),
\]
which proves \eqref{eq:variance} for finite sketch laws.

For a general sketch law \(\mu\), let
\(\mathcal S_1,\ldots,\mathcal S_m\) be independent draws from \(\mu\),
independent of \(U\). For each fixed \(U\),
\[
 \E_{\mathcal S_1,\ldots,\mathcal S_m}
 \left(\frac1m\sum_i\nolimits f(U\mathcal S_i)\right)
 =\int f(U\mathcal S)\,d\mu(\mathcal S).
\]
Since \(\E_U f(U\mathcal S_i)=\E_G f(\mathcal S_G)\)
for each fixed \(\mathcal S_i\), Jensen's inequality
and the preceding bound applied conditionally on
\(\mathcal S_1,\ldots,\mathcal S_m\) give
}
\begin{align*}
 \Var_U\!\left(\int f(U\mathcal{S})\,d\mu(\mathcal{S})\right)
 &\le \E_{\mathcal{S}_1,\ldots,\mathcal{S}_m}\Var_U\!\left(m^{-1}\sum_i\nolimits f(U\mathcal{S}_i)\right)\\
 &\le 2\,\E_{\mathcal{S}_1,\ldots,\mathcal{S}_m}
            \left\|m^{-1}\sum_i\nolimits P_{\mathcal{S}_i}\right\|_{\op}\Var_G(f(\mathcal S_G)).
\end{align*}
By the strong law of large numbers, for fixed $n$,
$
 \left\|m^{-1}\sum_i\nolimits P_{\mathcal S_i}\right\|_{\op}
 \longrightarrow \|M\|_{\op}=\rho_n
 \,\,\text{a.s. as }m\to\infty
$. 
Since these norms are bounded by one, dominated convergence
gives \eqref{eq:variance}.
\end{proof}

{
To extend \eqref{eq:variance} to \(V^T\), where $V$ is the Haar-distributed   basis in Introduction, let \(J\) be a diagonal
reflection with \(\det J=-1\), so that
\(O(n)=\SO(n)\cup J\SO(n)\). Since \(J\Sigma=\Sigma J\),
we have \(P_{\Sigma J\mathcal S}=JP_{\Sigma\mathcal S}J^T\) and hence
\(\mathcal R_k(\Sigma,J\mathcal S)=\mathcal R_k(\Sigma,\mathcal S)\).
Thus \(f(JU\mathcal S)=f(U\mathcal S)\) when
\(f=\mathcal R_k(\Sigma,\cdot)\) or any truncation thereof.
The Haar law of \(V^T\) is the equal mixture of the laws of \(U\) and
\(JU\). Consequently, for such \(f\),
\[
 \int f(V^T\mathcal S)\,d\mu(\mathcal S)
 \overset{d}{=}
 \int f(U\mathcal S)\,d\mu(\mathcal S).
\]
For each prescribed spectrum with \(\tau_k>0\),
\(\mathcal R_k(\Sigma,\cdot)\) is bounded and Borel measurable.
Taking \(f=\mathcal R_k(\Sigma,\cdot)\) in Lemma~\ref{lem:variance}
and using \eqref{eq:center} proves \eqref{eq:variance-comparison},
completing the proof of part (i).
\par}

\subsection{A moment bound at the scale of the optimal tail}

To obtain the bounds uniform over spectra in parts (ii) and (iii),
we need a uniform moment estimate for
\(\mathcal R_k(\Sigma,\cdot)\).
Bounding the numerator in \eqref{eq:residual} by
the total squared energy would introduce the unbounded ratio
\(\|\Sigma\|_F^2/\tau_k^2\). Instead we approximate the leading \(k\) rows
exactly and express the remaining error as a weighted sum over the tail.

\begin{lemma}\label{lem:moment}
For \(1<r\le2\) with \(2r<p+1\), \(Z\sim N(0,1)\), and \(Y\sim\chi^2_{p+1}\). Then
\begin{equation}\label{eq:moment}
  \|\mathcal R_k(\Sigma,\cdot)\|_{L^r(\nu)}
  \le K_{k,p,r} \coloneqq 1+k\bigl(\E|Z|^{2r}\,\E Y^{-r}\bigr)^{1/r}
  \qquad\text{for every }\Sigma\text{ with }\tau_k>0 .
\end{equation}
\end{lemma}

\begin{proof}
Let \(G\) be an \(n\times\ell\) standard Gaussian matrix,
so that \(\mathcal S_G=\range G\) has law \(\nu\).
Write \(B\) for the first \(k\) rows of \(G\) and \(z_i\)
for row \(i>k\). In the notation of Halko et al.
\cite[Theorem~9.1, Eq.~(9.3)]{halko2011}, take \(A=\Sigma\) and
\(\Omega=G\), so that \(\Omega_1=B\) and \(\Omega_2\) consists of the
rows \(z_i\), \(i>k\). Since \(B\) has full row rank almost surely,
the Frobenius-norm bound applies with \(\|\Sigma_2\|_F^2=\tau_k^2\) and
\(\|\Sigma_2\Omega_2\Omega_1^\dagger\|_F^2
=\sum_{i>k}\sigma_i^2\|z_iB^\dagger\|^2\).
Dividing by \(\tau_k^2\) gives
\begin{equation}\label{eq:weights}
 \mathcal R_k(\Sigma,\mathcal S_G)
 \le 1+\sum_{i>k} \nolimits w_i\|z_iB^\dagger\|^2,
 \qquad w_i \coloneqq  {\sigma_i^2}/{\tau_k^2},\qquad
 w_i\ge0,\quad \sum_{i>k}\nolimits w_i=1 .
\end{equation}
Note what division by \(\tau_k^2\) has done. Before normalization, the lower rows
carry the separate factors \(\sigma_i^2\). After normalization, their
contribution in \eqref{eq:weights} is a convex combination. In particular,
\begin{equation}\label{eq:convex-tail}
 \begin{aligned}
 \left\|\sum_{i>k}\nolimits w_i\|z_iB^\dagger\|^2\right\|_{L^r}
 &\le\sum_{i>k}\nolimits w_i\,\bigl(\E\|z_iB^\dagger\|^{2r}\bigr)^{1/r}
 =\bigl(\E\|z_{k+1}B^\dagger\|^{2r}\bigr)^{1/r}.
 \end{aligned}
\end{equation}
The inequality follows from Minkowski's inequality,
and the equality from the identical distributions of
\(z_iB^\dagger\) and \(\sum_{i>k}w_i=1\).
The resulting bound is independent of the spectrum.

It remains to bound one summand. Since \(z_i\) is a standard Gaussian
row independent of \(B\), conditional on \(B\), \(z_iB^\dagger\)
is a centered Gaussian row with covariance
\((B^\dagger)^TB^\dagger=(BB^T)^{-1}\). If \(\lambda_1,\ldots,\lambda_k\)
are the eigenvalues of this inverse, its squared norm has the law of
\(\sum_j\lambda_j Z_j^2\) for independent standard normals \(Z_j\).
Set \(s \coloneqq \sum_j\lambda_j\). By convexity with weights \(\lambda_j/s\),
\[
 \begin{aligned}
 \E\!\left[\|z_iB^\dagger\|^{2r}\mid B\right]
 &=\E_Z\left(\sum_j\nolimits\lambda_j Z_j^2\right)^r
 \le s^r\sum_j\nolimits(\lambda_j/s)\E|Z_j|^{2r}
 =\bigl(\tr[(BB^T)^{-1}]\bigr)^r\E|Z|^{2r}.
 \end{aligned}
\]
Taking expectation and applying \eqref{eq:iw} now gives
\[
  \bigl(\E\|z_iB^\dagger\|^{2r}\bigr)^{1/r}
         \le k\bigl(\E|Z|^{2r}\,\E Y^{-r}\bigr)^{1/r}.
\]
Inserting this bound into \eqref{eq:convex-tail} and adding the one
in \eqref{eq:weights} proves \eqref{eq:moment}.
\end{proof}

\subsection{Uniform \(L^1\) bounds}

We first establish the finite-dimensional bounds in part (ii).
To accommodate the cases where a uniform second-moment bound is
unavailable, we truncate the residual.

For a threshold \(a>0\), set
\(f(\mathcal S)\coloneqq\min\{\mathcal R_k(\Sigma,\mathcal S),a\}\).
Using
\(\mathcal R_k(\Sigma,\mathcal S)=f(\mathcal S)+
(\mathcal R_k(\Sigma,\mathcal S)-a)_+\),
\eqref{eq:center} and the triangle inequality give
\begin{align}
 \E_V|F_\Sigma(V)-g_\Sigma|
 &\le
 \underbrace{\E_V\left|\int f(V^T\mathcal{S})\,d\mu(\mathcal{S})-\int f\,d\nu\right|}
             _{\text{fluctuation of the truncated error}}
 +\underbrace{2\int
   (\mathcal R_k(\Sigma,\mathcal S)-a)_+\,d\nu(\mathcal S)}
             _{\text{discarded error}} .\label{eq:split}
\end{align}
Since \(V^T\mathcal S\) and \(\mathcal S_G\) have the same
law \(\nu\), the two truncation errors have equal expectations,
giving the factor two in \eqref{eq:split}.
{
Since \(\E_V\int f(V^T\mathcal S)\,d\mu=\int f\,d\nu\),
Cauchy--Schwarz and Lemma~\ref{lem:variance}, extended to \(V^T\) as above, give
\[
 \begin{aligned}
 \E_V\left|\int f(V^T\mathcal S)\,d\mu-\int f\,d\nu\right|
 &\le\sqrt{\Var_V\!\left(\int f(V^T\mathcal S)\,d\mu\right)}
 \le\sqrt{2\rho_n\,\E_G f(\mathcal S_G)^2}.
 \end{aligned}
\]
\par}
The pointwise bounds
\[
 \begin{aligned}
 f(\mathcal S)^2 \le a^{2-r}\mathcal R_k(\Sigma,\mathcal S)^r \,\,\,\,\text{and}\,\,\,\,
 (\mathcal R_k(\Sigma,\mathcal S)-a)_+
 \le a^{1-r}\mathcal R_k(\Sigma,\mathcal S)^r
 \end{aligned}
\]
give, by Lemmas~\ref{lem:variance} and~\ref{lem:moment},
\[
 \E_V |F_\Sigma(V)-g_\Sigma|
 \le (2\rho_n K_{k,p,r}^r a^{2-r})^{1/2}
          +2K_{k,p,r}^r a^{1-r}.
\]
For \(1<r<2\), the first term grows as \(a^{1-r/2}\),
whereas the second decays as \(a^{1-r}\).

Choosing \(a=K_{k,p,r}(2\rho_n)^{-1/r}\),
which is independent of \(\Sigma\), and taking the
supremum over spectra, we obtain
\begin{equation}\label{eq:explicit-rate}
 D_n\le 3K_{k,p,r}(2\rho_n)^{1-1/r}.
\end{equation}
This proves the general \(L^1\) bound, with
\(C_{k,p,r}=3K_{k,p,r}2^{1-1/r}\).

{
For \(p\ge4\), we can take \(r=2\) in Lemma~\ref{lem:moment}
and obtain a sharper constant without truncation.
By Cauchy--Schwarz, \eqref{eq:center}, \eqref{eq:variance-comparison},
and \eqref{eq:moment},
\[
 \E_V|F_\Sigma(V)-g_\Sigma|
 \le \sqrt{\Var_V F_\Sigma(V)}
 \le \sqrt{2\rho_n}\,\|\mathcal R_k(\Sigma,\cdot)\|_{L^2(\nu)}
 \le K_{k,p,2}\sqrt{2\rho_n}.
\]
Taking the supremum over spectra gives
\begin{equation}\label{eq:root-rate}
 D_n\le K_{k,p,2}\sqrt{2\rho_n},\qquad p\ge4.
\end{equation}
This is the stated square-root bound, with
\(C_{k,p}\coloneqq\sqrt{2}K_{k,p,2}\).  Since an admissible \(r>1\) exists for every \(p\ge2\),
\eqref{eq:explicit-rate} establishes the sufficient direction
of \eqref{eq:criterion}.
\par}

\subsection{One squared singular value detects a preferred direction}\label{sec:one-spike}

To prove the necessary direction of \eqref{eq:criterion},
consider a single-spike spectrum with squared singular values
\((n+1,1,\ldots,1)\). Then \(\tau_k^2=n-k\) and
\(K=A^TA=I+nvv^T\), where \(v=Ve_1\) is uniform on the
unit sphere. We keep the sketch width \(\ell=k+p\) fixed.

For a sketch space \(\mathcal S\) with orthonormal basis matrix
\(Q\), let \(P\coloneqq P_{\mathcal S}=QQ^T\),
\(z\coloneqq Q^Tv\), and \(t\coloneqq \|z\|^2=v^TPv\).
Since
\[
 Q^TKQ=I_\ell+nzz^T,\qquad
 Q^TK^2Q=I_\ell+n(n+2)zz^T,
\]
{
for \(z\ne0\), both compressed matrices act as the identity on \(z^\perp\).
Thus \eqref{eq:captured} gives
\[
 \|P_{A\mathcal S}A\|_F^2
 =\ell-1+\frac{1+n(n+2)t}{1+nt}
 =\ell+\frac{n(n+1)t}{1+nt}.
\]
The formula also holds for \(z=0\).
Using \(\tr K=2n\) and \(\tau_k^2=n-k\) in \eqref{eq:residual}, we obtain
}
\begin{equation}\label{eq:spike-loss}
 \mathcal R_k(A,\mathcal S)
 =\frac{n-\ell-1}{n-k}
 +\frac{n+1}{n-k}f(nt),\qquad
 f(x)\coloneqq\frac{1}{1+x}.
\end{equation}
Let \(J(v)\coloneqq\E_P f(nv^TPv)\). By
\eqref{eq:spike-loss},
$
 F_\Sigma(V)=\frac{n-\ell-1}{n-k}
 +\frac{n+1}{n-k}J(v)
$.
By \eqref{eq:center},
\(F_\Sigma(V)-g_\Sigma=\frac{n+1}{n-k}(J(v)-\E_vJ(v))\).
The definition of \(D_n\) therefore gives
\begin{equation}\label{eq:Jconcentrates}
 \E_v|J(v)-\E_vJ(v)|
 \le \frac{n-k}{n+1}D_n\le D_n.
\end{equation}

To relate the fluctuation of \(J(v)\) to \(\rho_n\),
let \(u\) be a unit eigenvector of \(M=\E P\)
corresponding to \(\rho_n\), and let
\(H(v)\coloneqq n(u^Tv)^2-1\).
We will express \(\operatorname{Cov}_v(J,H)\) as a multiple of
\(\rho_n-\ell/n\) and use it to bound \(D_n\) from below.

{
To compute the covariance, first fix \(P\).
By invariance of the conditional law of \(v\) given \(t\)
under orthogonal transformations preserving \(\range P\),
the conditional second moment has the form
\(\E_v[vv^T\mid t,P]=\alpha P+\beta(I-P)\).
Taking traces after multiplication by \(P\) and \(I-P\),
respectively, gives
\[
 \alpha\ell=\E_v[\|Pv\|^2\mid t,P]=t,\qquad
 \beta(n-\ell)=\E_v[\|(I-P)v\|^2\mid t,P]=1-t.
\]
Consequently,
\[
 \E_v[vv^T\mid t,P]=\frac{t}{\ell}P+
                             \frac{1-t}{n-\ell}(I-P).
\]
Let \(a\coloneqq u^TPu\) and \(X\coloneqq nt\).
}
Since \(H(v)=n(u^Tv)^2-1\) and \(X=nt\), we have
\[
 \begin{aligned}
 \E_v[H\mid X,P]
 &=\frac{aX}{\ell}+\frac{(1-a)(n-X)}{n-\ell}-1
 =\frac{n(a-\ell/n)}{\ell(n-\ell)}(X-\ell).
 \end{aligned}
\]
{
By \eqref{eq:sphere}, for every fixed \(P\),
\(X\) has the same law as
\(X_n\coloneqq n\sum_{j=1}^{\ell}v_j^2\),
with \(\E X_n=\ell\).
Since \(\E_v H=0\), we have
\begin{align}
 \operatorname{Cov}_v(J,H)
 &=\E_P\E_v[f(X)H(v)]
 =\E_P\E_v[f(X)\E_v(H\mid X,P)] \notag\\
 &=\frac{n}{\ell(n-\ell)}
 \E_P\left(a-\frac{\ell}{n}\right)
 \E[f(X_n)(X_n-\ell)] \notag\\
 &=\frac{n}{\ell(n-\ell)}
 \left(\rho_n-\frac{\ell}{n}\right)
 \operatorname{Cov}(f(X_n),X_n).
 \label{eq:spike-cov}
\end{align}
The third equality uses the independence of \(X\) and \(P\);
the last uses \(\E_P a=u^TMu=\rho_n\) and \(\E X_n=\ell\).
\par}

We next show that the covariance factor in \eqref{eq:spike-cov}
has a strictly negative limit.
Writing \(v=g/\|g\|\), we have
\(X_n=\sum_{j=1}^{\ell}g_j^2/(\|g\|^2/n)\).
Since \(\|g\|^2/n\to1\) in probability,
Slutsky's theorem gives \(X_n\overset{d}{\longrightarrow}
Y\sim\chi^2_\ell\).
Since \(xf(x)=1-f(x)\), \(\E X_n=\ell\),
and \(f\) is bounded and continuous,
\[
 \operatorname{Cov}(f(X_n),X_n)
 =1-(1+\ell)\E f(X_n)
 \longrightarrow 1-(1+\ell)\E f(Y)<0.
\]
The strict inequality follows from Jensen's inequality:
\(\E f(Y)>f(\E Y)=1/(1+\ell)\), since \(f\)
is strictly convex and \(Y\) is nonconstant.

{
To prove the lower bound on \(D_n\) in \eqref{eq:two-sided},
let \(b_\ell\coloneqq(1+\ell)\E f(Y)-1>0\).
By the preceding limit, \(|\operatorname{Cov}(f(X_n),X_n)|\ge b_\ell/2\)
for all sufficiently large \(n\). Since \(n/(n-\ell)\ge1\),
\eqref{eq:spike-cov} gives
\[
 |\operatorname{Cov}_v(J,H)|
 \ge\frac{b_\ell}{2\ell}\left(\rho_n-\frac{\ell}{n}\right).
\]
Since \(0\le J\le1\), we have \(\Var_v J\le\E_v|J-\E_v J|\).
By the sphere moments in Section~2.4,
\(\E_v H^2=2(n-1)/(n+2)\le2\).
Thus \eqref{eq:Jconcentrates} and Cauchy--Schwarz give
\[
 D_n\ge\E_v|J-\E_v J|\ge\Var_v J
 \ge\frac{|\operatorname{Cov}_v(J,H)|^2}{\E_v H^2}
 \ge\frac{b_\ell^2}{8\ell^2}\left(\rho_n-\frac{\ell}{n}\right)^2.
\]
This proves the lower bound in \eqref{eq:two-sided} with
\(c_\ell\coloneqq b_\ell^2/(8\ell^2)\).
In particular, \(D_n\to0\) implies \(\rho_n\to0\) for fixed \(\ell\),
completing the proof of part (ii).
\par}

\subsection{Uniform \(L^2\) bound}

{
The proof of part (iii) uses the same variance and moment
estimates as the square-root bound in Section~3.3.
For \(p\ge4\), by \eqref{eq:center},
\eqref{eq:variance-comparison}, and \eqref{eq:moment}
with \(r=2\),
\[
 \begin{aligned}
 \E_V|F_\Sigma(V)-g_\Sigma|^2
 &=\Var_V F_\Sigma(V)
 \le 2\rho_n\Var_G(\mathcal R_k(\Sigma,\mathcal S_G))\\
 &\le 2\rho_n\|\mathcal R_k(\Sigma,\cdot)\|_{L^2(\nu)}^2
 \le C_{k,p}^2\rho_n.
 \end{aligned}
\]
Taking the supremum over spectra and using
Cauchy--Schwarz, we obtain
\[
 D_n^2
 \le \sup_{\Sigma:\,\tau_k>0}\E_V|F_\Sigma(V)-g_\Sigma|^2
 \le C_{k,p}^2\rho_n.
\]
This proves \eqref{eq:uniform-second-moment}.
The equivalence follows from \eqref{eq:criterion},
completing the proof of Theorem~\ref{thm:criterion}.
\par}

{
\paragraph{Proof of Corollary~\ref{cor:non-haar}.} Recall that \(\nu_n\) is the Haar measure on \(O(n)\)
and \(\lambda_n\) is the law of \(V\).
By assumption, \(\lambda_n(B)\le C_n\nu_n(B)\)
for every Borel set \(B\subseteq O(n)\).
Consequently, for every nonnegative measurable function
\(f\) on \(O(n)\),
\[
 \int f\,d\lambda_n\le C_n\int f\,d\nu_n.
\]
Applying this inequality to the first and second powers of
$|F_\Sigma(V)-g_\Sigma|$ yields
\[
 \E_{V\sim\lambda_n}|F_\Sigma(V)-g_\Sigma|^q
 \le C_n\E_{V\sim\nu_n}|F_\Sigma(V)-g_\Sigma|^q,
 \qquad q=1,2 .
\]
For \(q=1\), taking the supremum over spectra and using
the definition of \(D_n\) in \eqref{eq:rho} gives
the first inequality in \eqref{eq:non-haar-bounds}.
For \(q=2\), by \eqref{eq:center} and
\eqref{eq:variance-comparison},
\[
 C_n\E_{V\sim\nu_n}|F_\Sigma(V)-g_\Sigma|^2
 =C_n\Var_{V\sim\nu_n}F_\Sigma(V)
 \le 2C_n\rho_n\Var_G\mathcal R_k(\Sigma,\mathcal S_G),
\]
which gives the second inequality in
\eqref{eq:non-haar-bounds}.
The sufficient direction and its convergence rates
then follow from Theorem~\ref{thm:criterion}
when \(C_n\) is bounded uniformly in \(n\).

If $c\nu_n\le\lambda_n\le C\nu_n$ for constants $0<c\le C<\infty$
independent of $n$, the two measure inequalities give
\[
 cD_n\le
 \sup_{\Sigma:\,\tau_k>0}\E_{V\sim\lambda_n}|F_\Sigma(V)-g_\Sigma|
 \le C D_n ,
\]
Thus convergence of the middle quantity to zero is equivalent to
$D_n\to0$, and
\eqref{eq:criterion} proves \eqref{eq:non-haar-criterion}.
}

\section{Consequences for common randomized sketches}\label{sec:applications}

In this section, we apply Theorem~\ref{thm:criterion}
to several commonly used randomized sketches.
We first identify two classes of sketches satisfying
\(\rho_n=\ell/n\), and then establish their agreement
on the same input. Numerical experiments illustrate
the theoretical results.

\subsection{Sampling rows of an orthogonal transform}

\begin{corollary}\label{cor:orthogonal-rows}
Let $O$ be an $n\times n$ orthogonal matrix, fixed or random. Given $O$,
choose a set $J$ of $\ell$ distinct rows so that each row has inclusion
probability $\ell/n$. Let $\mathcal{S}$ be the span of the chosen rows.
Then $\E P_{\mathcal{S}}=(\ell/n)I_n$. Consequently $D_n\to0$ as $n\to\infty$, with
$D_n\le K_{k,p,2}\sqrt{2\ell/n}$ when $p\ge4$.
\end{corollary}

\begin{proof}
Let \(R_J\) be the row selection matrix. Since
\(R_JOO^TR_J^T=I_\ell\), we have
\(P_{\mathcal S}=O^TR_J^TR_JO\).
Moreover, \(R_J^TR_J\) is diagonal with
\((R_J^TR_J)_{ii}=\mathbf{1}_{\{i\in J\}}\),
and \(\Pr(i\in J\mid O)=\ell/n\).
Consequently,
\begin{equation}\label{eq:uniform-rows}
 \E[P_{\mathcal S}\mid O]
 =O^T\E[R_J^TR_J\mid O]O
 =O^T\left(\frac{\ell}{n}I_n\right)O
 =\frac{\ell}{n}I_n.
\end{equation}
Taking expectation over \(O\) gives
\(\E P_{\mathcal S}=(\ell/n)I_n\), so
\(\rho_n=\ell/n\) by \eqref{eq:rho}. The convergence follows from
Theorem~\ref{thm:criterion}, and \eqref{eq:root-rate} gives the stated bound.
\end{proof}

{
Substituting $\rho_n=\ell/n$ into \eqref{eq:variance-comparison} gives
$
 \Var_V F_\Sigma(V)
 \le \frac{2\ell}{n}\Var_G\mathcal R_k(\Sigma,\mathcal S_G)
$.
For $p\ge4$, \eqref{eq:uniform-second-moment} and Chebyshev's inequality give
\[
 \sup_{\Sigma:\,\tau_k>0}
 \Pr_V\{|F_\Sigma(V)-g_\Sigma|>\varepsilon\}
 \le \frac{C_{k,p}^2\ell}{n\varepsilon^2}
 \qquad(\varepsilon>0).
\]
For \(p=2,3\), substituting \(\rho_n=\ell/n\)
into \eqref{eq:explicit-rate} gives
\(D_n=O(n^{-(1-1/r)})\).
Since \(1-1/r\to(p-1)/(p+1)\) as
\(r\to(p+1)/2\) from below, we obtain
\(D_n=O(n^{-1/3+\delta})\) and
\(D_n=O(n^{-1/2+\delta})\), respectively,
for every \(\delta>0\).
\par}

\paragraph{Examples.} Taking \(O=I\) gives uniform coordinate sampling. Taking \(O=HD\), where
\(H\) is the normalized Walsh--Hadamard transform and \(D\) contains
independent random signs on its diagonal, gives the real subsampled
randomized Hadamard transform (SRHT), in dimensions that are powers of
two \cite{tropp2011srht,boutsidis2013}. Taking \(O=CD_1\) or
\(O=CD_1CD_2\), with \(C\) the orthonormal type-II discrete cosine
transform (DCT) and \(D_1,D_2\) independent sign diagonals, gives one
and two signed cosine passes \cite{strang1999,epperlyThesis2025}.
Any further signed cosine passes also give an orthogonal matrix, so the
same calculation applies to their product.

\subsection{Sketches with sign and permutation symmetry}

For general sketch spaces, sign and permutation
symmetry also suffice to determine the mean projector \(\E P_{\mathcal S}\).

\begin{corollary}\label{cor:symmetries}
Let $\Omega$ be an $n\times\ell$ sketch of full column rank almost
surely. If the law of its column space $\mathcal{S} \coloneqq \range\Omega$ is invariant under coordinate
sign changes and coordinate permutations, then
$\E P_{\mathcal{S}}=(\ell/n)I_n$, and $D_n$ satisfies the same convergence and
rate as in Corollary~\ref{cor:orthogonal-rows}.
\end{corollary}

\begin{proof}
{
Let \(M\coloneqq\E P_{\mathcal S}\). Since
\(P_{T\mathcal S}=TP_{\mathcal S}T^T\) for every orthogonal \(T\),
the assumed symmetries give \(DMD^T=M\) and \(\Pi M\Pi^T=M\)
for every diagonal sign matrix \(D\) and permutation matrix \(\Pi\).
Choosing \(D\) to flip only coordinate \(i\) gives
\(M_{ij}=-M_{ij}=0\) for \(j\ne i\).
Choosing \(\Pi\) to interchange coordinates \(i\) and \(j\)
gives \(M_{ii}=M_{jj}\). Consequently,
\[
 M=\frac{\tr M}{n}I_n
 =\frac{\E\tr P_{\mathcal S}}{n}I_n
 =\frac{\ell}{n}I_n.
\]
Thus \(\rho_n=\ell/n\), and Theorem~\ref{thm:criterion}
and \eqref{eq:root-rate} give the stated conclusions.
\par}
\end{proof}

The variance comparison, probability bound, and
rates for $p=2,3$ following Corollary~\ref{cor:orthogonal-rows}
also apply to these sketches.

\paragraph{Examples.} Both symmetries hold for matrices with iid standard Gaussian entries or
iid entries uniform on $[-1,1]$, which have full column rank almost surely.
Matrices with independent uniform signs and the SparseIID and SparseStack
sketches used by Epperly and Webber also have these symmetries, although
they can lose rank. Since \(\operatorname{rank}(T\Omega)=\operatorname{rank}\Omega\)
for every signed permutation matrix \(T\), conditioning on
\(\operatorname{rank}\Omega=\ell\) preserves both symmetries whenever
this event has positive probability. Corollary~\ref{cor:symmetries}
then applies to the conditional law. Without conditioning,
\(\tr P_{\mathcal S}=\operatorname{rank}\Omega\), so the same argument gives
\(M=(\E\operatorname{rank}\Omega/n)I_n\).

{
\paragraph{Rank deficiency.} To see why rank deficiency matters, consider
\(A=\Sigma V^T\) with squared singular values
\(\sigma_i^2=\Gamma>1\) for \(i\le k\) and \(\sigma_i^2=1\) for \(i>k\).
If \(\operatorname{rank}\Omega<k\), then \(\dim(A\mathcal S)<k\),
so there exists a unit vector \(u\) in the leading left singular subspace
with \(u\perp A\mathcal S\). Hence
\[
 \|(I-P_{A\mathcal S})A\|_F^2
 \ge\|u^T(I-P_{A\mathcal S})A\|_2^2
 =\|u^TA\|_2^2=\Gamma.
\]
Dividing by \(\tau_k^2=n-k\) and averaging gives
\[
 \E_{\Omega,V}\mathcal R_k(A,\mathcal S)
 \ge\frac{\delta_n\Gamma}{n-k},
 \qquad \delta_n\coloneqq\Pr\{\operatorname{rank}\Omega<k\}.
\]
Let \(F_\Sigma(V)\) denote the conditional mean for the unconditioned
sketch as in \eqref{eq:means}. Then
\[
 \begin{aligned}
 \E_V|F_\Sigma(V)-g_\Sigma|
 &\ge\E_{\Omega,V}\mathcal R_k(A,\mathcal S)-g_\Sigma
 \ge\frac{\delta_n\Gamma}{n-k}-g_\Sigma.
 \end{aligned}
\]
By \eqref{eq:worst-case}, \(g_\Sigma\) is bounded independently of
\(\Gamma\). Letting \(\Gamma\to\infty\) gives
\(\sup_{\Sigma:\,\tau_k>0}\E_V|F_\Sigma(V)-g_\Sigma|=\infty\)
whenever \(\delta_n>0\).
Thus uniform convergence over spectra fails if \(\delta_n>0\)
for infinitely many \(n\), even if \(\delta_n\to0\).
The example for prescribed inputs in Section~\ref{sec:projection-limits}
shows a related effect even for a full-rank SRHT.
\par}

\paragraph{Remark (growing sketch width).}
{
The finite-dimensional bound \eqref{eq:root-rate} also allows
\(\ell\) to grow with \(n\) in Corollaries~\ref{cor:orthogonal-rows}
and~\ref{cor:symmetries}. For \(p=\ell-k\ge4\), we have \(\E Z^4=3\)
and \(\E Y^{-2}=1/[(p-1)(p-3)]\) by the inverse moment formula in
Section~\ref{sec:gaussian-tools}. Substituting these moments into
\eqref{eq:moment} with \(r=2\), and using \(\rho_n=\ell/n\)
in \eqref{eq:root-rate}, gives
\[
 D_n\le
 \left(1+k\sqrt{\frac{3}{(p-1)(p-3)}}\right)\sqrt{\frac{2\ell}{n}}
 \le(1+k)\sqrt{\frac{2\ell}{n}}.
\]
Thus \(D_n\to0\) for fixed \(k\), \(\ell=o(n)\), and \(p\ge4\)
eventually. The converse, however, depends on the width regime:
for a Gaussian sketch, \(D_n=0\) even when
\(\ell/n\to c\in(0,1)\) and hence \(\rho_n\to c\).
\par}

\subsection{Simultaneous mean error matching}

Suppose several sketches are available for the same matrix. Once we
know how close each conditional mean is to $g_\Sigma$, we can also
compare their mean errors on that matrix.

\begin{corollary}\label{cor:common-mean}
Fix $m\ge2$. For each $n$, let $\mathcal{S}^{(1)},\ldots,\mathcal{S}^{(m)}$ be random
$\ell$-dimensional sketch spaces with
$\rho_n^{(j)} \coloneqq \|\E P_{\mathcal{S}^{(j)}}\|_{\op}\to0$ for each $j=1,\ldots,m$ as
$n\to\infty$. For their conditional means evaluated at the same
Haar-distributed right singular basis $V$,
\[
 \sup_{\Sigma:\,\tau_k>0}\E_V
       \max_{1\le j\le m}|F_\Sigma^{(j)}(V)-g_\Sigma|\longrightarrow0 .
\]
In particular, for every pair $i,j$,
$\sup_{\Sigma:\,\tau_k>0}\E_V|F_\Sigma^{(i)}(V)-F_\Sigma^{(j)}(V)|\to0$.
\end{corollary}

\begin{proof}
Let \(D_n^{(j)}\) denote \(D_n\) in \eqref{eq:rho}
for sketch \(j\). By Theorem~\ref{thm:criterion},
\(D_n^{(j)}\to0\) for each \(j\). Since \(m\) is fixed,
\[
 \sup_{\Sigma:\,\tau_k>0}\E_V
 \max_{1\le j\le m}|F_\Sigma^{(j)}(V)-g_\Sigma|
 \le\sum_{j=1}^m D_n^{(j)}\longrightarrow0.
\]
Moreover, by the triangle inequality,
\[
 \sup_{\Sigma:\,\tau_k>0}\E_V
 |F_\Sigma^{(i)}(V)-F_\Sigma^{(j)}(V)|
 \le D_n^{(i)}+D_n^{(j)}\longrightarrow0.
\]
\end{proof}

\subsection{Numerical illustrations}

\subsubsection{Effect of the ambient dimension}

We first compare two block sampling schemes with different values of
\(\rho_n\) as the dimension \(n\) increases. Suppose \(\ell\) divides \(n\),
and partition the coordinate axes into \(n/\ell\) disjoint coordinate
spaces \(\mathcal S_1,\ldots,\mathcal S_{n/\ell}\), each of dimension
\(\ell\). We either select one block uniformly or always select the first
block. Their mean projectors are \((\ell/n)I_n\) and \(P_{\mathcal S_1}\),
so their respective values of \(\rho_n\) are \(\ell/n\) and \(1\).
For a fixed input \(A=\Sigma V^T\), both conditional means can be computed
exactly:
\[
 F_\Sigma^{(0)}(V)=\frac{\ell}{n}\sum_{b=1}^{n/\ell}
 \mathcal R_k(A,\mathcal S_b),\qquad
 F_\Sigma^{(1)}(V)=\mathcal R_k(A,\mathcal S_1).
\]
Under Haar rotations the two sketches have the same one-draw error law,
but their conditional means need not concentrate equally around \(g_\Sigma\).
By \eqref{eq:center} and \eqref{eq:uniform-second-moment}, for \(p\ge4\),
\[
 \sqrt{\Var_V F_\Sigma(V)}
 =\big(\E_V|F_\Sigma(V)-g_\Sigma|^2\big)^{1/2}
 \le C_{k,p}\sqrt{\rho_n}.
\]

For this experiment and the next, we use \(k=4\), \(\ell=8\), and three
spectra: a single spike with \(\sigma_1=\sqrt{n+1}\) and \(\sigma_i=1\)
for \(i>1\); a polynomial spectrum with \(\sigma_i=1/i\); and a
shrinking-tail spectrum with \(\sigma_i=1/i\) for \(i\le k\) and
\(\sigma_{k+j}=1/(nj)\) for \(1\le j\le n-k\).
The last spectrum has \(\tau_k^2=\Theta(n^{-2})\) and tests the matching
of relative errors on the scale of a vanishing optimal tail.

\begin{figure}[tb]
 \centering
 \includegraphics[width=\linewidth]{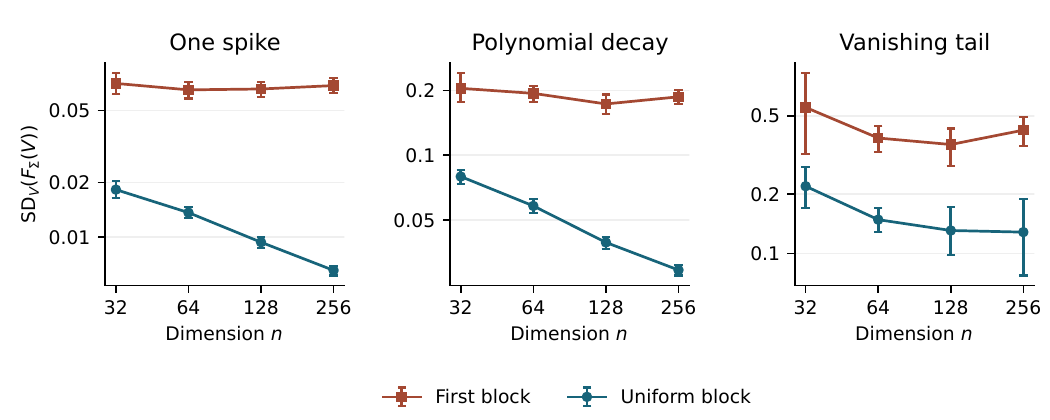}
 \caption{Variation of the conditional mean across 600 independent Haar draws of the right
 singular basis at each dimension. For each basis and spectrum, we compute
 both the first-block error and the mean over all disjoint blocks exactly.
 Points show sample standard deviations, and bars give 95\% percentile
 bootstrap intervals obtained by resampling the bases. Both axes are logarithmic.}\label{fig:blocks}
\end{figure}

Figure~\ref{fig:blocks} shows that, for the single spike, increasing
\(n\) from 32 to 256 leaves the first-block standard deviation near
\(0.07\), whereas the uniform-block standard deviation decreases
from \(0.0183\) to \(0.00655\). The polynomial spectrum shows a similar
separation. The shrinking-tail spectrum has the widest intervals,
with the last two uniform-block estimates close to each other.
The results illustrate the effect of \(\rho_n\) on concentration, while
the magnitude of the fluctuations also depends on the spectrum.

\subsubsection{Effect of the mean projector norm}

We next vary \(\rho_n\) at a fixed dimension by mixing the two preceding
sketch laws. With probability \(\alpha\in[0,1]\) select the first block,
and otherwise select a block uniformly. The mean projector is
\(\alpha P_{\mathcal S_1}+(1-\alpha)(\ell/n)I_n\).
Hence its largest eigenvalue and the conditional mean of the mixture are
\begin{equation}\label{eq:block-mixture}
 \rho_n=\alpha+(1-\alpha)\frac{\ell}{n},\qquad
 F_\Sigma^{(\alpha)}(V)=(1-\alpha)F_\Sigma^{(0)}(V)+\alpha F_\Sigma^{(1)}(V).
\end{equation}
Thus both quantities can be evaluated exactly for every \(\alpha\).

\begin{figure}[tb]
\centering
\includegraphics[width=\linewidth]{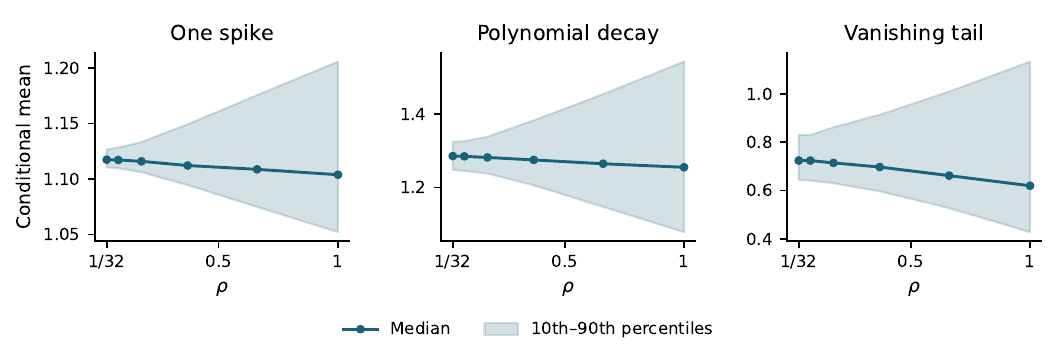}
\caption{Conditional means as the sketch increasingly favors one block.
Each point is the median over the same 600 Haar draws of the right singular basis
at $n=256$, $k=4$, $\ell=8$. The band runs from the 10th to the 90th
percentile. For each basis, \eqref{eq:block-mixture} gives the
conditional mean exactly. The vertical scales differ across panels.}
\label{fig:mixture}
\end{figure}

Figure~\ref{fig:mixture} shows how the central 80\% interval of
\(F_\Sigma^{(\alpha)}(V)\) changes as \(\alpha\) increases from zero
to one. For the single spike, the interval widens from
\([1.110,1.127]\) to \([1.052,1.206]\). At a given \(\alpha\),
the three spectra have the same \(\rho_n\) but different spreads.
Thus \(\rho_n\) controls concentration uniformly over spectra without
determining the variance for a particular spectrum.

\subsubsection{Variance comparison at fixed \(\rho_n\)}

We now compare uniform coordinate sampling, signed DCT, and SRHT.
By Corollary~\ref{cor:orthogonal-rows}, all three have
\(\rho_n=\ell/n\). The question is whether their actual variances are
comparable at a fixed dimension.

Unlike the block schemes above, exact evaluation of \(F_\Sigma(V)\)
for these sketches is costly. To separate variation across \(V\) from
the noise of a finite sketch average, for each \(V\) we generate two
conditionally independent estimates \(Y_1,Y_2\) satisfying
\(\E[Y_1\mid V]=\E[Y_2\mid V]=F_\Sigma(V)\).
Since \(\operatorname{Cov}(Y_1,Y_2\mid V)=0\), the law of total covariance
gives
\[
 \operatorname{Cov}(Y_1,Y_2)
 =\operatorname{Cov}\big(\E[Y_1\mid V],\E[Y_2\mid V]\big)
 =\Var_V F_\Sigma(V).
\]
We can therefore estimate the variance of the conditional mean without
including the conditional sampling variance of either replicate.

\begin{table}[tb]
\centering
\caption{Estimated variance of the conditional mean for the single spike at
\(n=256\), \(k=4\), \(\ell=8\). All entries are in units of \(10^{-6}\).
The covariance is estimated over 1,024 independent draws of the spike direction. Each of
its two replicates averages over 128 independent random partitions into
blocks, averaging the losses within each partition. A block in a random
partition is a uniform sample of eight coordinates. The signed transforms
also redraw their signs for each partition. Intervals are 95\% percentile
bootstrap intervals for the paired observations.}\label{tab:transforms}
\smallskip
\begin{tabular}{lrr}
\toprule
Sketch & $\Var_V(F_\Sigma(V))$ & Interval\\
\midrule
Uniform coordinates & \(9.68\) & \([8.80,10.69]\)\\
Signed DCT          & \(0.015\) & \([-0.006,0.036]\)\\
SRHT                & \(0.009\) & \([-0.011,0.032]\)\\
\bottomrule
\end{tabular}
\end{table}

In Table~\ref{tab:transforms}, the confidence interval for coordinate
sampling lies above zero, whereas those for both signed transforms
include zero and have upper endpoints far below the coordinate estimate.
The signed-transform variances are thus much smaller in this experiment,
although the estimates do not resolve their precise magnitudes.
This shows that sketches with the same \(\rho_n\) can exhibit very
different fluctuations at a given dimension.

\subsubsection{Comparison of common randomized sketches}

Finally, we examine the six sketches and two spectra of Epperly
and Webber's Figure~2 under our Haar input model
\cite[Sections~2.2 and~3.2]{epperlyWebber2026}.
They varied the sketch width at a fixed dimension; here we vary \(n\)
with \(k=10\) and \(\ell=24\) fixed.
For SparseIID, each vector \(\Omega^Te_i\) has 16 nonzeros on average;
for SparseStack, it has one signed entry in each of eight disjoint
blocks \cite[Definitions~1.6 and~3.2]{camanoOSI}.
At each dimension, we use the same Haar bases for all six sketches
and compare their estimated conditional means with the spectrum-specific
Gaussian mean, rather than the worst-case bound in \eqref{eq:worst-case}.

\begin{figure}[!htbp]
\centering
\includegraphics[width=\linewidth]{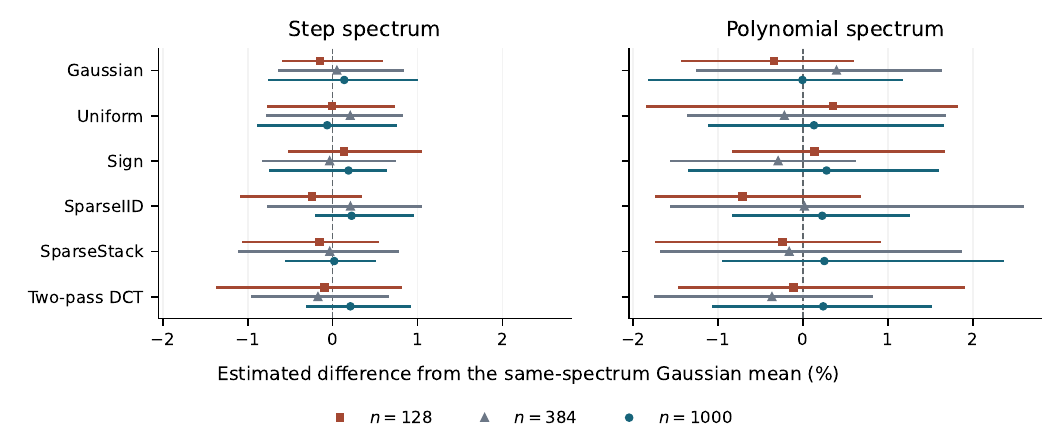}
\caption{The six sketches on the two spectra of Epperly and Webber,
with a Haar right singular basis. The step spectrum has ten unit values
and a flat $10^{-5}$ tail. The polynomial spectrum has $\sigma_i=i^{-2}$.
At each dimension, all six sketches use the same 24 independent
bases. Each conditional mean uses 128 fresh draws with $k=10$ and
$\ell=24$. Points and bars give the median and the 10th and 90th
percentiles of $100(\widehat F_\Sigma(V)/\widehat g_\Sigma-1)$ over
those bases. The Gaussian reference $\widehat g_\Sigma$ uses 8,192
separate draws for each dimension and spectrum.}
\label{fig:ew-haar}
\end{figure}

At \(n=1000\), all six medians in Figure~\ref{fig:ew-haar} are within
\(0.3\%\) of the Gaussian reference for both spectra. For the two-pass
DCT, the central 80\% intervals are \([-0.30,0.91]\%\) for the step
spectrum and \([-1.05,1.50]\%\) for the polynomial spectrum, comparable
to those of the Gaussian sketch. These intervals include sketch
sampling noise: even for a Gaussian sketch, whose true conditional mean
is \(g_\Sigma\) for every \(V\), a finite sample average fluctuates.
Unlike the paired covariance estimates in Table~\ref{tab:transforms},
the intervals in Figure~\ref{fig:ew-haar} do not isolate the variation
due to the choice of basis. The discrete sketches were sampled
unconditionally, and all draws in this experiment had full column rank.

\section{Discussion}\label{sec:discussion}

We first illustrate the non-Haar extension in
Corollary~\ref{cor:non-haar} and examine its limitations. We then
distinguish Gaussian limits for individual sketch projections from
Gaussian mean error matching, both under the Haar input model and
for prescribed singular directions.

\subsection{Beyond Haar right singular bases}\label{sec:non-haar}

{
\paragraph{An example for sufficiency.} Corollary~\ref{cor:non-haar} allows input distributions with a genuine
directional preference. For example, fix $0<a<1$ and draw
\[
 x\sim N(0,\diag(a,1,\ldots,1)).
\]
Take $x/\|x\|$ as the first right singular vector and choose a uniform
orthonormal completion in its orthogonal complement. Compared with the
standard Gaussian density, the density of $x$ has ratio
\[
 a^{-1/2}\exp\bigl(-\tfrac12(a^{-1}-1)x_1^2\bigr)
 \le a^{-1/2}.
\]
The same bound holds for the density of the resulting basis relative
to Haar measure: normalization and random completion use the same
conditional rule for both Gaussian distributions, and the standard
Gaussian construction produces a Haar basis. Thus the sufficient
direction and its rates in Corollary~\ref{cor:non-haar} apply, with
a density bound independent of dimension. Yet
\[
 \sqrt n\,V_{11}\ \Longrightarrow\ N(0,a)
\]
because $\|x\|^2/n\to1$ in probability. Under Haar measure this limit
would be $N(0,1)$. The difference between the two limiting laws
also shows that these basis distributions do not approach Haar in
total variation.

\paragraph{An example for the equivalence.}
 To obtain a non-Haar distribution satisfying both density bounds of
Corollary~\ref{cor:non-haar}, use a Haar basis with probability
$\theta\in(0,1)$ and the preceding basis with probability $1-\theta$,
where $\theta$ is fixed. Its density relative to Haar lies between
$\theta$ and $\theta+(1-\theta)a^{-1/2}$. The limit law of
$\sqrt n\,V_{11}$ is now the corresponding mixture of $N(0,1)$ and
$N(0,a)$, still different from the standard normal law. This gives
a class of non-Haar inputs for which the full geometric
characterization and the convergence rates remain valid.

Haar measure conditioned on an event of uniformly positive Haar
probability also satisfies the upper density bound.

\paragraph{Insufficiency of total variation convergence.} Conversely, convergence to Haar in total variation alone does not
ensure mean error matching uniformly over spectra. To see this, use uniform
coordinate sampling with $k=1$, and draw the input basis from Haar
measure with probability $1-\delta_n$ and set it to $I_n$ with
probability $\delta_n$, where $\delta_n>0$ and $\delta_n\to0$.
The total variation distance to Haar is at most $\delta_n$, but the
mean absolute difference is at least
$\delta_n|F_\Sigma(I_n)-g_\Sigma|$. For every fixed $n>\ell$, the
single-spike spectrum in Section~1.2 makes this lower bound
arbitrarily large as the gap grows. Thus even a vanishing probability
of a special input can prevent agreement uniformly over spectra.
\par}

\subsection{Can a central limit theorem guarantee mean error matching?}\label{sec:projection-limits}

{
One might wonder whether a central limit theorem for sketch
projections alone suffices to establish mean error matching.
The following two counterexamples show that it does not.

\paragraph{Haar inputs.}
Consider \(\Sigma=\diag(\sqrt{n+1},1,\ldots,1)\) and \(v=Ve_1\),
where \(V\) is Haar distributed.
Let \(\Omega^{(1)}\) select the first \(\ell\) coordinates and
\(\Omega^{(2)}\) select \(\ell\) distinct coordinates uniformly,
independently of \(V\). Then \(\rho_n^{(1)}=1\) and
\(\rho_n^{(2)}=\ell/n\).
By permutation invariance and the representation \(v=g/\|g\|\),
\[
 \sqrt n\,(\Omega^{(1)})^Tv
 \overset{d}{=}\sqrt n\,(\Omega^{(2)})^Tv
 \overset{d}{\longrightarrow}N(0,I_\ell).
\]
For the fixed-coordinate sketch, \eqref{eq:spike-loss} applies to its
conditional mean \(F_\Sigma^{(1)}(V)\) with
\(nt=n\|(\Omega^{(1)})^Tv\|^2\overset{d}{\longrightarrow}Y\sim\chi_\ell^2\).
By \eqref{eq:center} and boundedness of this single-spike error,
\(g_\Sigma\to1+\E(1+Y)^{-1}\).
For the uniform-coordinate sketch, Theorem~\ref{thm:criterion} gives
\(\E_V|F_\Sigma^{(2)}(V)-g_\Sigma|\to0\).
Consequently,
\[
 \begin{aligned}
 F_\Sigma^{(1)}(V)&\overset{d}{\longrightarrow}1+\frac{1}{1+Y},
 &\quad&\text{fixed coordinates},\\
 F_\Sigma^{(2)}(V)&\overset{L^1}{\longrightarrow}1+\E\frac{1}{1+Y},
 &&\text{uniform coordinates}.
 \end{aligned}
\]
The first limit is nonconstant, whereas the second is the Gaussian mean
limit. Thus even identical one-draw projection laws do not determine
the concentration of the conditional mean.

\paragraph{Prescribed inputs.}
For a prescribed input, a further obstruction is that a Gaussian
projection limit need not ensure convergence of the mean error.
Consider the SRHT \(\Omega^T=R_JHD\) from
Section~\ref{sec:applications}, with \(n=2^m\), \(k=1\), fixed
\(\ell>2\), and \(v=n^{-1/2}(1,\ldots,1)^T\) as the leading
right singular vector. For any fixed row set \(J\) and unit
\(a\in\R^\ell\),
\[
 a^T\sqrt n\,\Omega^Tv=\sum_{j=1}^n c_jD_{jj},
 \qquad c_j=(a^TR_JH)_j,
 \qquad \sum_{j=1}^n c_j^2=1,\quad
 \max_j|c_j|\le\sqrt{\ell/n}.
\]
The Lindeberg condition holds for every sequence of row sets.
By Cram\'er--Wold,
\[
 \sqrt n\,\Omega^Tv\overset{d}{\longrightarrow}N(0,I_\ell).
\]
The limit also holds for uniformly sampled rows independent of \(D\).
For related sketch central limit theorems, see Ahfock, Astle, and
Richardson \cite[Theorem~2]{ahfock2021}.

Now choose a single-spike spectrum with squared singular values
\((\Gamma_n,1,\ldots,1)\), where \(\Gamma_n=n^2 2^n\).
On the event \(D=I_n\) and \(1\notin J\), we have
\(\Omega^Tv=R_JHv=R_Je_1=0\). This event has probability
\(2^{-n}(1-\ell/n)\).
Since \(A^TAv=\Gamma_nv\), an exact miss gives \(Av\perp A\mathcal S\)
and hence \(\|(I-P_{A\mathcal S})A\|_F^2\ge\|Av\|^2=\Gamma_n\).
With \(\tau_1^2=n-1\), it follows that
\[
 \E_{\mathcal S}\mathcal R_1(A,\mathcal S)
 \ge 2^{-n}(1-\ell/n)\frac{\Gamma_n}{n-1}
 =\frac{(1-\ell/n)n^2}{n-1}\longrightarrow\infty.
\]
The Gaussian mean remains bounded by Lemma~\ref{lem:moment}.
Thus the projection CLT alone does not ensure uniform integrability
of the relative errors on prescribed inputs: rare misses can dominate
their expectation.
\par}

\section*{Use of AI}
Experiments in randomized SVD have often found that quite different sketches
give errors comparable to those of Gaussian sketches
\cite{halko2011,epperlyWebber2026}.
Like many researchers, the author wondered what might explain this agreement.
Several months ago, the author explored the question with ChatGPT-5.6 Sol, but
those discussions yielded only partial results for specific sketches.

On November, the author returned to the question after reading  \cite{openaiLift}.
To compare the spectral prices of two congruent configurations, the proof
maps both into a higher-dimensional space with the same Gaussian map and
uses its polar factor to recover an isometry. The author saw a possible
parallel with the way a single input rotation acts on all the sketch spaces
and chated with ChatGPT-6.2 Sol whether a similar lift could help explain the
mean error. Those discussions led to a result for two-level spectra,   as in the \emph{step} example of  
\cite[Figure~2]{epperlyWebber2026}. 
The author then asked ChatGPT how to handle an arbitrary spectrum and it suggested using the Brascamp--Lieb inequality  
to bound the variance of the mean error under a shared rotation.
Over several rounds of discussion and revision with ChatGPT and Codex,
 the argument for sufficiency and the single-spike
argument for necessity were developed, arriving at the geometric criterion proved here.
ChatGPT   also assisted with drafting, numerical illustrations, and 
  Lean formalization of   the results, which is available at  
\url{https://github.com/yuningyang19/rsvd-typical-orientations-lean}.

\appendix
\section{Unbounded worst-case mean error of signed DCT sketches}\label{app:paired-dct}

Fix \(1\le k<\ell<n\). Suppose
\(\Pr(\mathcal S=\mathcal S_0)=\eta>0\) for some sketch space \(\mathcal S_0\).
Choose a unit vector \(u\perp\mathcal S_0\) and a \(k\)-dimensional
subspace \(\mathcal V\) containing \(u\).
For \(\Gamma>1\), let \(A=\sqrt{\Gamma}P_{\mathcal V}+I-P_{\mathcal V}\).
Since \(Au=\sqrt{\Gamma}u\), we have
\(\langle u,As\rangle=\langle Au,s\rangle=0\) for every \(s\in\mathcal S_0\).
Thus \(u\perp A\mathcal S_0\), which gives
\[
 \|(I-P_{A\mathcal S_0})A\|_F^2
 \ge\|(I-P_{A\mathcal S_0})Au\|^2
 =\|Au\|^2=\Gamma.
\]
As \(\tau_k^2=n-k\), it follows that
\[
 \E_{\mathcal S}\mathcal R_k(A,\mathcal S)
 \ge\frac{\eta\Gamma}{n-k}\longrightarrow\infty
 \qquad(\Gamma\to\infty).
\]
The two-pass DCT sketch has finite support at each fixed \(n\), so such
a space \(\mathcal S_0\) exists. This proves the second equality in
\eqref{eq:worst-case}.

\end{document}